\documentclass[a4paper,11pt]{article}

\usepackage{graphicx}
\usepackage{bm}
\usepackage[percent]{overpic}

\usepackage{authblk}
\usepackage{amscd}
\usepackage{amsfonts}
\usepackage[numbers]{natbib}
\usepackage{amsthm}
\usepackage{amssymb}
\usepackage{amsmath}
\usepackage{algpseudocode}
\usepackage{algorithm}
\usepackage{cancel}
\usepackage{chngcntr}
\usepackage{caption}
\usepackage{enumerate}
\usepackage{epsfig}
\usepackage{epic}
\usepackage{eepic}
\usepackage{enumitem}
\usepackage{etoolbox}
\usepackage{fancyhdr}
\usepackage{fullpage}
\usepackage[a4paper, total={6in, 8in}]{geometry}
\usepackage{graphicx}
\usepackage{here}
\usepackage[utf8]{inputenc}
\usepackage{lscape}
\usepackage{listings}
\usepackage{longtable}
\usepackage{natbib}
\usepackage{nccmath}
\usepackage{parskip}
\usepackage{rotating}
\usepackage{subfigure}
\usepackage{tabularx}
\usepackage{textcomp}
\usepackage{tikz}
\usepackage{titlesec}
\usepackage{threeparttable}
\usepackage{varwidth}
\usepackage{framed,multirow}
\graphicspath{ {./Figures/} }

\usepackage{xcolor}

\title{Asymmetry and condition number of an elliptic-parabolic system for biological network formation}

\author[1]{Clarissa Astuto}
\author[1,2]{Daniele Boffi}
\author[1]{Jan Haskovec}
\author[1,3]{Peter Markowich}
\author[4]{Giovanni Russo}

\affil[1]{King Abdullah University of Science and Technology (KAUST), 4700, Thuwal, Saudi Arabia}
\affil[2]{Department of Mathematics ''F. Casorati'', University of Pavia, Pavia, Italy}
\affil[3]{Department of Mathematics, University of Vienna, Vienna, Austria}
\affil[4]{Department of Mathematics and Computer Science, University of Catania, Catania, Italy}

\begin{document}
\maketitle

\begin{abstract}
We present results of numerical simulations of the tensor-valued elliptic-parabolic PDE model for biological network formation.
The numerical method is based on a non-linear finite difference scheme on a uniform Cartesian grid in a 2D domain. 
The focus is on the impact of different discretization methods and choices of regularization parameters
on the symmetry of the numerical solution.
In particular, we show that using the symmetric alternating-direction implicit (ADI) method for time discretization helps preserve the symmetry of the solution, compared to the (non symmetric) ADI method.
Moreover, we study the effect of regularization by isotropic background permeability $r>0$, showing that increased condition number of the elliptic problem due to decreasing value of $r$ leads to loss of symmetry. {We show that in this case, neither the use of the symmetric ADI method preserves the symmetry of the solution.}
Finally, we perform numerical error analysis of our method {making use of} Wasserstein distance.
\end{abstract}

\section{Introduction}Principles of formation, adaptation and functioning of 
biological transportation networks have been a long standing topic of scientific investigation {for} significant applications in leaf venation in plants {malinowski2013understanding},
vascular pattern formation \cite{sedmera2011function}, mammalian circulatory systems or
neural networks that transport electric charge \cite{eichmann2005guidance,michel1995morphogenesis}.
Typical subjects of investigation are geometrical and topological properties of optimal networks, their statistical properties and robustness with respect to damage or varying external conditions.
For instance, in mammalian circulatory systems one aim of study is the relation between the dilation of arteries and an augmentation of blood flow \cite{pohl1986crucial}. Other studies reveal that local gradient of pressure can interfere with the diameter of blood vessels as an adaptive response to the stress \cite{hacking1996shear,pries1998structural,chen2012haemodynamics,hu2012blood}.

In plant leafs, the pattern of their venation seems to influence the cells that are engaged in photosynthesis, and other functionalities of the plant, such as its longevity and the optimal water distribution.
Modeling of formation and adaptation of leaf venation is a very challenging task because of the nature of the problem. Every leaf of the same plant exhibits
different venation patterns \cite{hu2013adaptation}.
This is reflected by the inherent non-uniqueness of solutions and, even, instabilities in the corresponding mathematical models.
Consequently, it is difficult to validate results of numerical simulations versus experimental observations.
A small change in the parameters of the model or its discretization (such as the resolution of the numerical grid) can lead to very different solutions.

The modeling framework for biological network formation 
introduced by Hu and Cai in \cite{hu2013adaptation, hu2013optimization} involves a purely local dynamic adaptation model based on mechanical laws, consisting of a system of ordinary differential equations (ODE) on graph edges coupled to a linear system of equations for the material pressure.
The biological nature of the model is reflected by a metabolic cost function that is proportional to a power of the conductance of the edge. 
Local conservation of mass is imposed by the Kirchhoff law. The model responds to merely local information and naturally incorporates fluctuations in flow distributions. In \cite{hu2013optimization} a related PDE-based continuum model was proposed, which consists of a parabolic reaction-diffusion equation for the vector-valued network conductivity, coupled to a Poisson equation for the pressure. The model was subsequently studied in the series of papers \cite{marko_perthame,marko_perthame_2, marko_perthame_schlo, marko_albi, albi_burger,portaro}.
A more general model with tensor-valued conductivity was proposed in {\cite{hu2019optimization} and further studied in \cite{marko_pilli}}. The {vector- and tensor-valued} modeling approaches were compared numerically in \cite{astuto2022comparison}.

This paper focuses on numerical treatment of the tensor-valued PDE model of \cite{marko_pilli}, where the permeability tensor appearing in the Poisson equation is regularized by adding a multiple $r\mathbb{I}$ of the identity matrix, with $r>0$.
We discretize the system in space using a finite-difference scheme on a two-dimensional Cartesian grid. The time discretization is carried out using the alternating-direction implicit (ADI) and symmetric-ADI schemes. The time discretization is crucial since the  system is stiff in all its components. {In system (\ref{eq_darcy_p_fin}-\ref{eq_reaction_diff_fin}) we will see a diffusion term $\Delta \mathbb{C}$, that needs to be treated implicitly, a nonlinear term $||\mathbb{C}||^{\gamma - 2}\mathbb{C}$, with $\gamma < 1$ and the variable $\mathbb{C}$ that assumes values close to 0, and least but not last, there is a pressure term $\nabla p \otimes \nabla p$, whose eigenvalues are 0 and $|\nabla p|$, and it can generate stiffness when $|\nabla p|$ is very large.}

The stiffness may lead to loss of symmetry of the solution in situations when all parameters, initial and boundary data are symmetric. The main goal of the paper is to investigate the loss of symmetry and its dependence on the model parameters, in particular on the value of the regularization parameter $r$.
Moreover, we carry out convergence analysis of the method. Here we argue that the Wasserstein distance~\cite{carrillo,otto1996double,ottoF} is an appropriate choice of distance in the convergence analysis. {The numerical solution takes very small values in a significant portion of the computational domain, away from the branches of the emergent network. Since the shape, position and number of the branches change with the discretization (number of points of the discrete domain), using an $L^p$ norm is not advantageous in studying the convergence analysis. Instead, we propose to use the Wasserstein distance for performing the convergence analysis, as it addresses this problem by considering the displacement of the solutions obtained with different numerical grids.}

The paper is organized as follows: In Section~\ref{section_model} we introduce the tensor-valued PDE model. In Section~\ref{section_methods} we describe the finite-difference semi-implicit schemes that we use for its discretization.
In Section~\ref{section_results} we provide the results of the numerical simulations, with focus on investigation of the symmetry of the solution and its dependence on the time discretization and the parameter values. Finally, in Section~\ref{section_conclusion} we summarize the results and draw some conclusions.

\section{The PDE Model}
\label{section_model}
The PDE model, {proposed in \cite{hu2019optimization} and further studied in \cite{marko_pilli},} consists of an elliptic equation for the pressure $p = p(t,\vec{x}) \in \mathbb{R}$ representing the Darcy's law, and a parabolic reaction-diffusion equation for the tensor-valued conductivity $\mathbb{C} = \mathbb{C}(t,\vec{x})$,
\begin{align}
	\label{eq_darcy_p_fin}
	-\nabla\cdot \left( (r\mathbb{I}+\mathbb{C} {)}\nabla p \right) &= S, \\
	\label{eq_reaction_diff_fin}
	\frac{\partial\mathbb{C}}{\partial t} - D^2\Delta \mathbb{C} -c^2\nabla p \otimes \nabla p + \alpha ||\mathbb{C}||^{\gamma-2} \mathbb{C} & = 0.
\end{align}
The term $S = S(\vec{x})$ denotes the distribution of sources and sinks, which has to be prescribed as a datum. The function $r :\Omega \to \mathbb{R}^+$, with $r(x) \geq {r_0} > 0$, describes the isotropic background permeability of the medium.
In \eqref{eq_reaction_diff_fin} the diffusion coefficient $D>0$ controls the random effects in the transportation medium and the activation parameter $c^2 > 0$ describes the tendency of the network to align with the pressure gradient.
The reaction term $\alpha ||\mathbb{C}||^{\gamma-2} \mathbb{C}$ models the metabolic cost of maintaining the network structure, with metabolic coefficient $\alpha>0$, metabolic exponent $\gamma>0$, {and where $||\mathbb{C}||$ denotes the Frobenius norm of the matrix $\mathbb{C}$.}
For blood circulatory systems we choose $\gamma=1/2$, see \cite{hu2012blood} for details, while for modeling of leaf venation in plants we have $1/2\leq \gamma \leq 1$, see \cite{hu2013adaptation, hu2013optimization}.

We pose \eqref{eq_darcy_p_fin}--\eqref{eq_reaction_diff_fin}
on a bounded domain $\Omega \subset \mathbb{R}^2$ with smooth boundary $\partial \Omega$.
We choose homogeneous Neumann boundary conditions for $\mathbb{C}$ and $p$ on $\partial\Omega$,
\begin{equation}
	\label{eq_bc}
	\nabla \mathbb{C}(t,\vec{x})\cdot \nu = 0, \quad {\nu\cdot \left( (r\mathbb{I}+\mathbb{C} )\nabla p \right) = 0}, \quad \vec{x}\in \partial \Omega, \, t\geq 0
\end{equation}
where $\nu$ is the outer normal vector to $\partial \Omega$ and the boundary condition for $\mathbb{C}$ is interpreted elementwise, i.e., 
$\nabla \mathbb{C}_{i,j}(t,\vec{x})\cdot \nu=0$
for all $i,j=1,\dots,d$.
{The choice of Neumann boundary conditions for $\mathbb{C}$ is motivated by the fact that in \cite{astuto2022comparison} we observed that homogeneous Dirichlet boundary conditions for $\mathbb{C}$ lead to the formation of boundary layers in the solutions. The choice of homogeneous Neumann boundary conditions suppresses this undesirable effect; we further comment on this issue in Section \ref{section_space_discrete}. Moreover, let us note that the same type of boundary condition was chosen in the numerical simulations carried out in \cite{hu2019optimization}.
	We also note that the reaction-diffusion equation \eqref{eq_reaction_diff_fin} is, after an eventual regularization of the metabolic term $\alpha ||\mathbb{C}||^{\gamma-2} \mathbb{C}$ if $\gamma\leq 1$, see \eqref{eq_Qcal}, is well-posed  subject to \eqref{eq_bc}.
}

{The choice of homogeneous Neumann boundary conditions for the pressure $p$ is natural from the modeling point of view, as it imposes zero flux of the transported material through the boundary $\partial\Omega$. We note that the material flux is given by $\vec q = (r\mathbb{I}+\mathbb{C} )\nabla p$.}
Moreover, we observe that the homogeneous Neumann boundary condition imposes the global mass balance
\begin{equation}
	\label{eq_Source}
	\int_\Omega S(\vec{x}) d\vec{x} = 0
\end{equation}
as a necessary condition to ensure solvability of \eqref{eq_darcy_p_fin}.
{Indeed, integrating \eqref{eq_darcy_p_fin} over $\Omega$ and using the Green formula and \eqref{eq_bc} gives
	\begin{equation*}
		\int_\Omega S(\vec{x}) d\vec{x} =
		\int_{\partial\Omega}
		\nu\cdot \left( (r\mathbb{I}+\mathbb{C} )\nabla p \right) d\vec{x} = 0.
	\end{equation*}
	Existence of a solution $p\in H^1(\Omega)$, unique up to an additive constant, is then a direct consequence of the Lax-Milgram lemma.
}

Finally, we prescribe a positive semidefinite initial condition $\mathbb{C}^0 \geq 0$ for the conductivity $\mathbb{C}$,
\begin{equation}
	\label{eq_ic}
	\mathbb{C}(t=0,\vec{x}) = \mathbb{C}^0(\vec{x}) \quad \text{ in } \Omega.
\end{equation}
A fundamental observation about the system \eqref{eq_darcy_p_fin}--\eqref{eq_reaction_diff_fin}
is that it represents an $L^2$-gradient flow of the energy
\begin{equation}  \label{Energy}
	\mathcal{E}[\mathbb{C}] = \int_\Omega \frac{D^2}{2} |\nabla\mathbb{C}|^2 + c^2\nabla p[\mathbb{C}] \cdot (r\mathbb{I} + \mathbb{C}) \nabla p[\mathbb{C}] + \frac{\alpha}{\gamma} ||\mathbb{C}||^\gamma \, d\vec{x},
\end{equation}
where $p[\mathbb{C}]$ is the unique (up to an additive constant) solution of \eqref{eq_darcy_p_fin} subject to the homogeneous Neumann boundary condition \eqref{eq_bc}.
{
	Proof of this claim is obtained by a simple modification of \cite[Lemma 1]{marko_perthame_2}.
	Indeed, calculation of the Fr\'echet derivative
	of the first and third terms in \eqref{Energy},
	i.e., the diffusive and metabolic terms, is straightforward.
	For the convenience of the reader, we demonstrate here the calculation of the first-order Fr\'echet derivative of the kinetic term  with respect to $\mathbb{C}$ in direction $\Phi\in [H^1(\Omega)]^{d\times d}$,
	\begin{equation}  \label{GF_calc}
		\frac{\delta}{\delta\mathbb{C}} \int_\Omega \nabla p[\mathbb{C}] \cdot (r\mathbb{I} + \mathbb{C}) \nabla p[\mathbb{C}] \, d\vec{x} =
		2 \int_\Omega \nabla p[\mathbb{C}] \cdot (r\mathbb{I} + \mathbb{C}) \nabla \frac{\delta p[\mathbb{C},\Phi]}{\delta\mathbb{C}} \, d\vec{x}
		+ \int_\Omega \nabla p[\mathbb{C}] \cdot \Phi \nabla p[\mathbb{C}] \, d\vec{x}.
	\end{equation}
	The first-order variation of the weak formulation of the Poisson equation \eqref{eq_darcy_p_fin} in direction $\Phi \in [H^1(\Omega)]^{d\times d}$
	with test function $\phi\in H^1(\Omega)$ reads
	\begin{equation*}
		\int_\Omega \nabla\phi\cdot (r\mathbb{I} + \mathbb{C}) \nabla \frac{\delta p[\mathbb{C},\Phi]}{\delta\mathbb{C}} + \nabla\phi\cdot \Phi \nabla p[\mathbb{C}] \, d\vec{x} = 0,
	\end{equation*}
	and choosing $\phi:=p[\mathbb{C}]$ yields
	\begin{equation*}
		\int_\Omega \nabla p[\mathbb{C}] \cdot (r\mathbb{I} + \mathbb{C}) \nabla \frac{\delta p[\mathbb{C},\Phi]}{\delta\mathbb{C}} \, d\vec{x} =
		-    \int_\Omega \nabla p[\mathbb{C}] \cdot \Phi \nabla p[\mathbb{C}] \, d\vec{x}.
	\end{equation*}
	Using this in \eqref{GF_calc}, we finally obtain
	\begin{align*}
		\frac{\delta}{\delta\mathbb{C}} \int_\Omega \nabla p[\mathbb{C}] \cdot (r\mathbb{I} + \mathbb{C}) \nabla p[\mathbb{C}] \, d\vec{x} &=
		- \int_\Omega \nabla p[\mathbb{C}] \cdot \Phi \nabla p[\mathbb{C}] \, d\vec{x} \\
		&=
		- \int_\Omega \left( \nabla p[\mathbb{C}] \otimes \nabla p[\mathbb{C}] \right) : \Phi\, d\vec{x},
	\end{align*}
	where $\mathbb{A}:\mathbb{B}$ denotes the scalar product of the matrices $\mathbb{A}$, $\mathbb{B}$.
}

In \cite{marko_pilli} it has been shown that for $\gamma>1$ the energy functional \eqref{Energy} is coercive and strictly convex. Consequently, it possesses a unique minimizer that describes the optimal transportation structure for the given distribution of sources and sinks $S$. On the other hand, {and how we said at the beginning of this section, in the specific application of leaf venation the parameter $\gamma$ belongs to the interval $(0,1)$} and it renders the energy highly non-convex with a multitude of critical points. {The main difficulty is the negative exponent $\gamma - 2$ in the reaction term, because in practise we are dividing by 0 when the conductivity variable $\mathbb{C}$ assumes negligible values. For this reason we add a '\textit{stabilization parameter}' in Eq.~\eqref{eq_Qcal}.} This fact is manifested in the numerical simulations carried out in this paper, where we shall observe their strong sensitivity with respect to the choice of the initial datum $\mathbb{C}^0$.

\section{Numerical schemes}
\label{section_methods}
In this section we briefly describe a fully second order space and time discretization that we adopt in our numerical simulations. We refer to \cite{astuto2022comparison} for further details.

{We adopt a finite differences scheme in space and a semi-implicit scheme in time, where, for the time-advancing, we use the Alternating Direction Implicit method, in symmetric form. The considered equations have been already analyzed, from the numerical point of view, making use of finite elements methods (FEM), in \cite{albi_burger} for the vectorial form (where the unknown conductivity variable is $m\in \mathbb{R}^2$ and where the authors only consider the case of $r = 0.1$), and in \cite{astuto2023finite} for the tensor model. For the tensor model, a FEM monolithic scheme, in mixed formulation with adaptive mesh refinement is already under investigation.}


\subsection{Space discretization}
\label{section_space_discrete}
For the discretization in space we choose the two-dimensional quadratic domain $\Omega = [0,1]\times [0,1]$ where we construct a uniform Cartesian mesh with spatial step $h := \Delta x = \Delta y$. We denote $\Omega_h$ the discrete computational domain. The discretized conductivity, ${ \mathbb{C}_{i,j}\approx \mathbb{C}(x_i,y_j)}$, and pressure, ${ p_{i,j}\approx p(x_i,y_j)}$, are defined at the center of the cell $(i,j)$, therefore we have $x_i=(i-1/2)h, \, y_j=(j-1/2)h, \,(i,j)\in\{1,\dots,N\}^2$, $h N = 1$. {With the choice of a cell centered discretization, it becomes natural to impose the homogeneous Neumann boundary conditions, with the technique of the eliminated ghost-points, and it guarantees the exact conservation of the total mass of the solution. Here we show a simple 1D example to show how it works. A generic diffusion equation can be written as $\partial_t u = -\partial_x J$, where the expression for the flux is $J = - \partial_x u$ and $x \in [0,1]$. If we impose homogeneous Neumann boundary conditions in $x = 0$, in the proximity of the boundary, we have
	\[ 
	\partial_t u_1 = - (\partial_x J)_1 = - \frac{J_{3/2} - J_{1/2}}{h}.
	\]
	where $u_i \approx u(x_i)$ and $J_i = - (u_{i+1/2} - u_{i-1/2})/h$. Since we have $J_{1/2} = 0$ from the boundary conditions, we simply consider $\partial_t u_1 = - J_{3/2}/h$.} 

{To justify the choice of the boundary conditions, in \cite{astuto2022comparison} we show an anomalous behavior of the solution near to the boundaries, in the case of homogeneous Dirichlet boundary conditions and zero-diffusivity (i.e., $D = 0$). We also show that \textit{ad-hoc} boundary conditions can be derived from the model, with the complications of adding nonlinearity to the equations. To overcome these difficulties, here we impose homogeneous Neumann boundary conditions, as it has been done in \cite{hu2019optimization}.}


The space discretization of Eq.~\eqref{eq_reaction_diff_fin} written in compact form (see \cite{astuto2022comparison} for more details) reads
\begin{eqnarray}
	\label{eq_C_compact}	
	\displaystyle \frac{\partial {\mathbb{C}}_{\rm comp}}{\partial t}&=& D^2\mathcal{L}\,{\mathbb{C}}_{\rm comp} + c^2 \mathcal{P} - \alpha \mathcal{Q}(\mathbb{C}) {\mathbb{C}}_{\rm comp}
\end{eqnarray}
where $\mathcal{L}$ is the discrete, {second order}, Laplacian operator,
${\mathbb{C}}_{\rm comp} = [C^{(1,1)},C^{(1,2)},C^{(2,2)}]^T$ is the vector of the unknowns for the conductivity, and $\mathcal{P}$ is the
matrix of pressure gradients $\mathcal{P} = [ \mathcal{D}_x p\,\mathcal{D}_x p, \mathcal{D}_x p\,\mathcal{D}_y p, \mathcal{D}_y p\,\mathcal{D}_y p ]$, {where $\mathcal{D}_x$ and, resp., $\mathcal{D}_y$ are
	the discrete first-order derivative operators in
	the $x-$ and, resp., $y-$direction, with central difference approximation, that ensures the second order accuracy.}  For the metabolic terms, we have
\begin{eqnarray}
	\label{eq_Qcal}
	\mathcal{Q}(\mathbb{C}) &=& ||\mathbb{C} { + \varepsilon} ||^{\gamma - 2}.
\end{eqnarray}
where $||\cdot||$ denotes the Frobenius norm { and $\varepsilon > 0$ is a stabilization parameter to avoid the division by 0 when $\mathbb{C}$ is close to 0. In \cite{astuto2022comparison} we justify the choice of the value of this parameter in our tests, i.e., $\varepsilon = 10^{-3}$.}

The {Poisson equation \eqref{eq_darcy_p_fin} can be extended} as
\begin{align}
	\partial_x\left(\left(r + C^{(1,1)} \right) \partial_x p \right) + \partial_x\left( C^{(1,2)} \partial_y p \right) + \partial_y\left( C^{(1,2)} \partial_x p \right)  
	+ \partial_y\left(\left(r + C^{(2,2)} \right) \partial_y p \right) = -S,
\end{align}
{while the discrete version is
	\begin{align}
		\mathcal{D}_x\left(\left(r + C^{(1,1)} \right) \mathcal{D}_x p \right) + \mathcal{D}_x\left( C^{(1,2)} \mathcal{D}_y p \right) + \mathcal{D}_y\left( C^{(1,2)} \mathcal{D}_x p \right)  
		+ \mathcal{D}_y\left(\left(r + C^{(2,2)} \right) \mathcal{D}_y p \right) = -S,
\end{align}}
where we use the symmetry of the conductance tensor $C^{(1,2)} = C^{(2,1)}$. 
We discretize the components of the above formula one by one, since we use different discretizations for each term. For simplicity of notation we define $\mathcal{C}^{(1,1)} = r + C^{(1,1)}$ and $\mathcal{C}^{(2,2)} = r + C^{(2,2)}$. Then we have
\begin{align} \nonumber    {\mathcal{D}}_x\left(\mathcal{C}^{(1,1)} {\mathcal{D}}_x\, p \right)_{i,j} \approx \frac{1}{2h^2} \left( \left( \mathcal{C}^{(1,1)}_{i+1,j} + \mathcal{C}^{(1,1)}_{i,j}\right)p_{i+1,j} + \left( \mathcal{C}^{(1,1)}_{i-1,j} + \mathcal{C}^{(1,1)}_{i,j}\right)p_{i-1,j} - \right. \\ \left. \left( \mathcal{C}^{(1,1)}_{i+1,j} + \mathcal{C}^{(1,1)}_{i-1,j} + 2\mathcal{C}^{(1,1)}_{i,j}\right)p_{i,j}  \right)  
\end{align}
We omit the term with both $y-$derivatives because it is analogous to the one with $x-$derivatives.  The term with mixed derivatives is discretized as follows,
\begin{align*}
	&{\mathcal{D}}_x\left(C^{(1,2)} {\mathcal{D}}_y\, p \right)_{i,j} \approx  \frac{1}{8 h^2}\left( C^{(1,2)}_{i+1,j} + C^{(1,2)}_{i,j}\right)(p_{i+1,j+1} - p_{i+1,j-1})  
	\\ &-\frac{1}{8 h^2} \left(  \left( C^{(1,2)}_{i-1,j} + C^{(1,2)}_{i,j}\right)(p_{i-1,j+1} - p_{i-1,j-1}) + \left( C^{(1,2)}_{i+1,j} - C^{(1,2)}_{i-1,j} \right) (p_{i,j+1} - p_{i,j-1}) \right)
\end{align*}
and analogously for the term with the $y,x$-derivatives.

\subsection{Time discretization: \textit{symmetric}-ADI method and extrapolation technique}
In this section we describe the {semi-implicit} time discretization that we apply to the model. It is a crucial point since the Eq.~\eqref{eq_reaction_diff_fin} is very stiff in all its components, {as we said in the Introduction.}  


We choose the \textit{symmetric} alternating-direction implicit (ADI) scheme {(see, for instance, \cite{alcubierre1992time}), and we treat explicitly the nonlinearity in the metabolic term.} 
From its definition, the classical ADI scheme \cite{peaceman_rachford} is not symmetric, since we choose which direction considering implicit for the first half step, while the second one is automatically chosen. {In this way, it allows us to solve linear systems of dimension $N\times N$ for each direction, instead of dimension $N^2\times N^2$, drastically reducing the computational cost of computing the solution of the equations for the three components of the tensor $\mathbb{C}$. } A symmetrized version of this method will compute the average of the two choices. 

{Since the Eq.~\eqref{eq_darcy_p_fin} for the pressure depends also on the conductivity $\mathbb{C}$, the second order accuracy is not guaranteed for the solution of the system (\ref{eq_darcy_p_fin}-\ref{eq_reaction_diff_fin}). For this reason we consider an improvement of the ADI scheme, that we already described in \cite{CiCP-31-707}.} We extrapolate the solution of the reaction-diffusion equation to compute the solution of the Poisson equation for the pressure: given the conductivity tensor at time $t^n$ and $t^{n-1}$, we extrapolate the conductivity at time $t^{n+1/2}$
\begin{eqnarray}
	\label{eq_extrapol}	\displaystyle \mathbb{C}^{n+1/2} = \frac{3}{2} \mathbb{C}^n - \frac{1}{2} \mathbb{C}^{n-1},
\end{eqnarray}
then we compute $p^{n+1/2}$ by 
solving the Poisson equation $-\nabla\left( \left( \mathbb{C}^{n+1/2} \right) \nabla p^{n+1/2} \right)  = S,$ which in the discretized version reads
\begin{equation}
	\label{eq_poisson_discrete}
	-\mathcal{L}\left(\mathbb{C}^{n+1/2}\right)\, p^{n+1/2} = S.
\end{equation}
Finally, we apply the \textit{symmetric}-ADI method to solve~\eqref{eq_C_compact}. {In practise we apply twice the traditional ADI scheme. The first time we start with the $y-$direction implicit and the $x-$direction explicit, and it reads
	\begin{eqnarray*}
		1^{\rm st}-{\rm step} \qquad	\left( I - \frac{\Delta t}{2}\mathcal{L}_y \right) \widetilde{\mathbb{C}}_1 &=&  \left(I + \frac{\Delta t}{2}\mathcal{L}_x \right) \mathbb{C}^n  + \Delta t\,\mathcal{P}^{n} \\
		2^{\rm nd}-{\rm step} \qquad  	\left( I - \frac{\Delta t}{2}\mathcal{L}_x + \Delta t\,\alpha \mathcal{Q}^c\left({ \mathbb{C}^{n}}\right)\right) \mathbb{C}^{n+1}_y&=& \left( I + \frac{\Delta t}{2}\mathcal{L}_y \right) \widetilde{\mathbb{C}}_1
	\end{eqnarray*} 
	where $ \displaystyle \mathcal{L}_\beta $, with $  \beta=x,y$, { are the discrete operators for the second derivatives} in $x$ and $y$ directions, resp., with $  \displaystyle \mathcal{L}_\beta \in  \mathbb{R}^{N\times N}$. In the $1^{\rm st}-$step we solve for $\widetilde{\mathbb{C}}_1$,\footnote{{The numerical scheme is implemented in Matlab and the solution of the linear system is computed with the '\textbackslash' command.}} and in the $2^{\rm nd}$ one we solve for $\mathbb{C}^{n+1}_y$, and so far we have the traditional ADI scheme. Now, with the opposite order, we start considering the $x-$direction implicit, and the $y-$direction explicit, and, analogously in the $2^{\rm nd}-$step we solve for $\mathbb{C}^{n+1}_x$.
	At the end we calculate $\mathbb{C}^{n+1}$ as the average between the two solutions, s.t., $\mathbb{C}^{n+1} = \frac{1}{2}{ \mathbb{C}^{n+1}_x} + \frac{1}{2}{ \mathbb{C}^{n+1}_y}.$ }

A numerical comparison  between the classical ADI and the symmetric ADI methods is presented in Table~\ref{table_symmetry_ADI2}.

\subsubsection{Comparisons between ADI and \textit{symmetric}-ADI method}
In this section we compare the two versions of the ADI method to see the improvements in the symmetric version. We adopt a symmetric numerical scheme, and we choose symmetric initial datum and source function $S$.
Consequently, the exact solution to the problem retains symmetry at each time step.

In order to check if our scheme is symmetric, we calculate the asymmetry of the solution with the following formula
\begin{equation}
	\label{eq_asymm}
	{\rm asymm}(A) = \frac{||A - A^T||}{||A + A^T||}
\end{equation}

In Table~\ref{table_symmetry_ADI2} we see the comparison between the two ADI schemes, for different choices of the regularization parameter $r$, representing the background permeability in the elliptic equation \eqref{eq_darcy_p_fin}.

\begin{table}[h]
	\centering 
	\caption{\textit{In this table we see the difference between a traditional ADI scheme and the symmetric version, for different values of the background permeability $r$. The choice of the parameters, the source function $S$ and the initial datum are defined in Eqs.~(\ref{eq_ic1_def}-\ref{eq_choice_param}), {with $\vec{x}_0 = (0.25,0.25)$}.}}
	\begin{tabular}{ccccccc}
		\hline
		& \multicolumn{2}{c}{$r = 10^{-2}$} & \multicolumn{2}{c}{$r = 10^{-3}$} & \multicolumn{2}{c}{$r = 10^{-4}$} \\
		\hline         & ADI & sym-ADI & ADI & sym-ADI & ADI & sym-ADI  \\ \hline
		asymm($\mathbb{C}$)  & 6.107e-04 & 4.390e-08 & 2.591e-02 &  1.673e-03 & 1.536e-01 &  1.868e-01 \ \\ \hline
		asymm($p$)  & 3.8862e-05  &  4.289e-09 &  6.137e-02 & 2.688e-03 &  2.683e-02 & 2.808e-02  \\ \hline
	\end{tabular}
	\label{table_symmetry_ADI2}
\end{table}  


Moreover, investigating the reasons why we lose the symmetry of the solution, we notice that it is related to the computation of the solution of the Poisson equation \eqref{eq_poisson_discrete} and in the choice of the background permeability $r$. When this parameters tends to zero, the condition number of the iteration matrix $\mathcal{L}(\mathbb{C})$ increases, up to the order $10^{8}$. In Fig.~\ref{fig_condition} we show the quantity defined in Eq.~\eqref{eq_asymm} for the module of the conductivity tensor $\mathbb{C}$ and the pressure $p$, as functions of time, together with the condition number of the matrix $\mathcal{L}(\mathbb{C})$. {We show how the asymmetry of the solutions increases when we decrease the parameters $r$.}



\section{Numerical Results}
\label{section_results}
The numerical results focus on the effect of the regularization parameter $r$ in \eqref{eq_darcy_p_fin}, and how some properties of the numerical scheme are strictly connected to it.


In our simulations we define the following initial conditions and source function $S$, as in equations \eqref{eq_Source} and \eqref{eq_ic}:
\begin{align}
	\label{eq_ic1_def}
	C^{(1,1)}(t=0) = 1, \quad C^{(1,2)}(t=0) = 0 \quad C^{(2,2)}(t=0) = C^{(1,1)}\\
	\label{eq_ic3_def}
	C^{(1,1)}(t=0) = (2 - |x + y|)\exp(-10|x-y|),\quad C^{(1,2)}(t=0) = 0, \quad C^{(2,2)}(t=0) = C^{(1,1)} \\ \label{eq_ic2_def}
	C^{(1,1)}(t=0) = 0, \quad C^{(1,2)}(t=0) = 0 \quad C^{(2,2)}(t=0) = C^{(1,1)} \\
	S(\vec{x}) = E - \bar{E}, \quad  E = \exp(-\sigma(\vec{x}-\vec{x}_0)^2), \quad \sigma = 500, \quad \vec{x}_0 = (0.25,0.25)
\end{align}
where $\bar{E} = {\rm mean}(E)$. The values of parameters of the system that we used in the simulations are as follows:
\begin{equation}
	\label{eq_choice_param}
	\alpha = 0.75, \, c = 5, \, D = 10^{-2}, \, \varepsilon = 10^{-3},
\end{equation}
with number of points $N = 600$, time step $\Delta t = h$ and final time $t_{\rm fin} = 10$. The choice of the parameter $\varepsilon$ is justified in the recent paper \cite{astuto2022comparison}, where we show that the solutions are qualitatively very close for $\varepsilon = 10^{-3}$ and $\varepsilon = 10^{-4}$.

\begin{figure}[tb]
	\centering
	\begin{minipage}{.32\textwidth}  \begin{overpic}[abs,width=1.1\textwidth,unit=1mm,scale=.25]{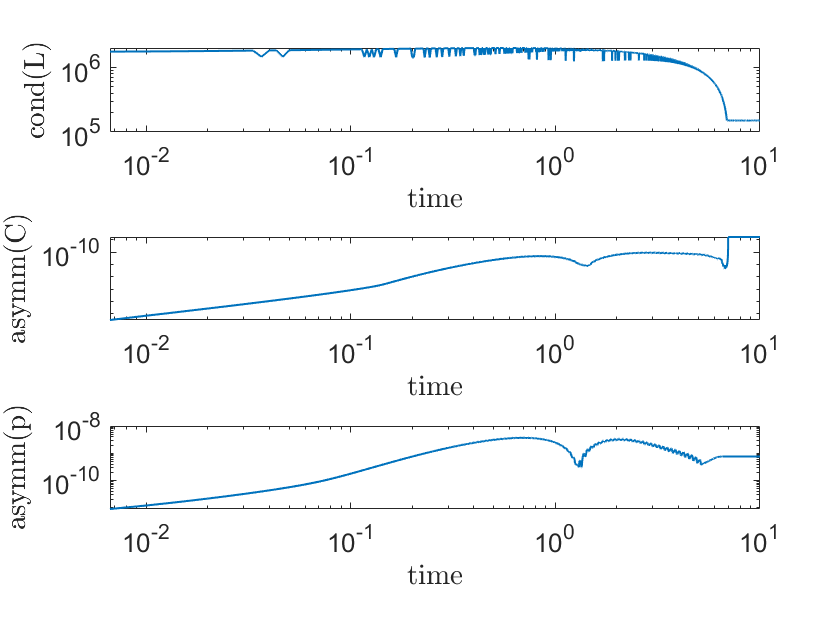}
			\put(12.,39){N = 600, $r = 10^{-1}$}
		\end{overpic}	  
	\end{minipage}
	\begin{minipage}{.32\textwidth}  
		\begin{overpic}[abs,width=1.1\textwidth,unit=1mm,scale=.25]{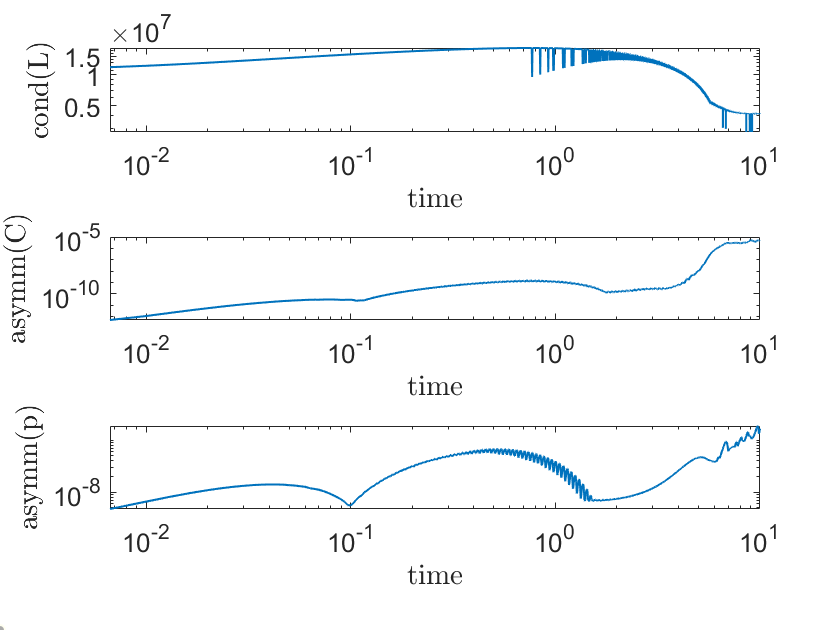}
			\put(12.,39){N = 600, $r = 10^{-2}$}
		\end{overpic}    
	\end{minipage}
	\begin{minipage}{.32\textwidth}    \begin{overpic}[abs,width=1.1\textwidth,unit=1mm,scale=.25]{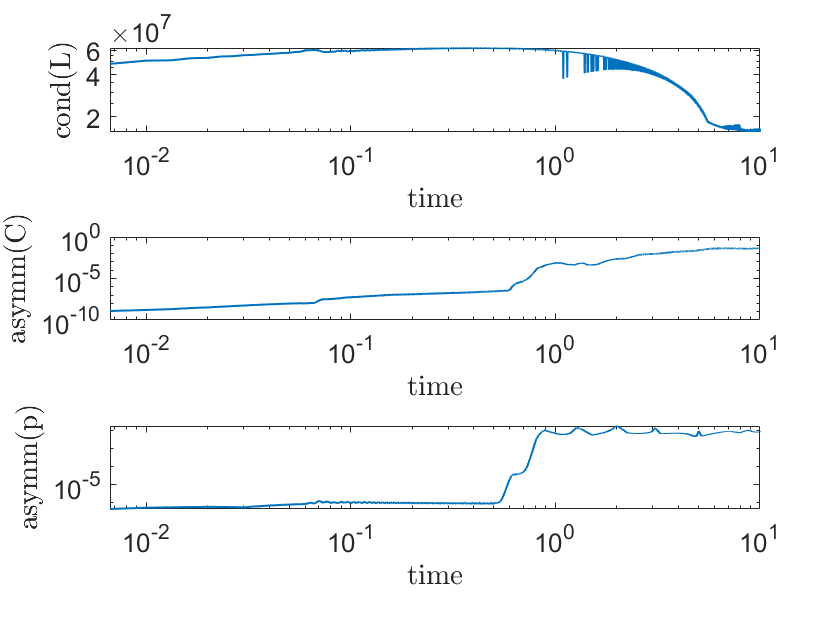}
			\put(12.,39){N = 600, $r = 10^{-3}$}
		\end{overpic}    
	\end{minipage}
	\caption{\textit{Quantity defined in {Eq.~\eqref{eq_asymm}} for the module of the conductivity tensor $\mathbb{C}$ and the pressure $p$, as function of time, together with the condition number of $\mathcal{L}$, for different values of $r = 10^{-1},10^{-2},10^{-3}$. { In this test $\vec{x}_0 = (0.1,0.1)$.}}}
	\label{fig_condition}
\end{figure}

In Fig.~\ref{fig_condition} {and in Table~\ref{table_symmetry_ADI2}} we show how the symmetry of the solutions strongly depends on the parameter $r$, {for two different initial conditions, $\vec x_0 = (0.25,0.25)$ in Fig.~\ref{fig_condition} and $\vec x_0 = (0.1,0.1)$ in Table~\ref{table_symmetry_ADI2}}. We calculate the asymmetry of the two variables conductivity $\mathbb{C}$ and pressure $p$, defined in Eq.~\eqref{eq_asymm}, at each time step, and we compare these quantities with the conditioning number of the elliptic operator for the Poisson equation $\mathcal{L}(\mathbb{C})$. Since the Eq.~\eqref{eq_darcy_p_fin} is strongly degenerate for $r\to 0$, it is not possible to consider negligible values for the parameter. {In Fig.~\ref{fig_condition} we see that for $r = 10^{-1}$ the symmetry of the solutions is well guaranteed, with the asymmetry of the order of $10^{-10}$ for the conductivity and $10^{-8}$ for the pressure. For $r = 10^{-2}$ the asymmetry is of the order of $10^{-5}$ for the conductivity and $10^{-8}$ for the pressure, while for $r = 10^{-3}$ we completely lose the symmetry for the conductivity, and, respectively, we see the growing of the values for the condition number of the Laplacian operator (see cond(L) in the same plot). Analogously, we see the asymmetry of the two solutions growing in Table~\ref{table_symmetry_ADI2} for $r = 10^{-2}, 10^{-3}, 10^{-4}$. For $r = 10^{-4}$ the symmetry is completely lost, losing also the advantages of the symmetric version of the ADI scheme.}
\begin{figure}[H]
	\centering
	\begin{minipage}{.32\textwidth}    \includegraphics[width=1.1\textwidth]{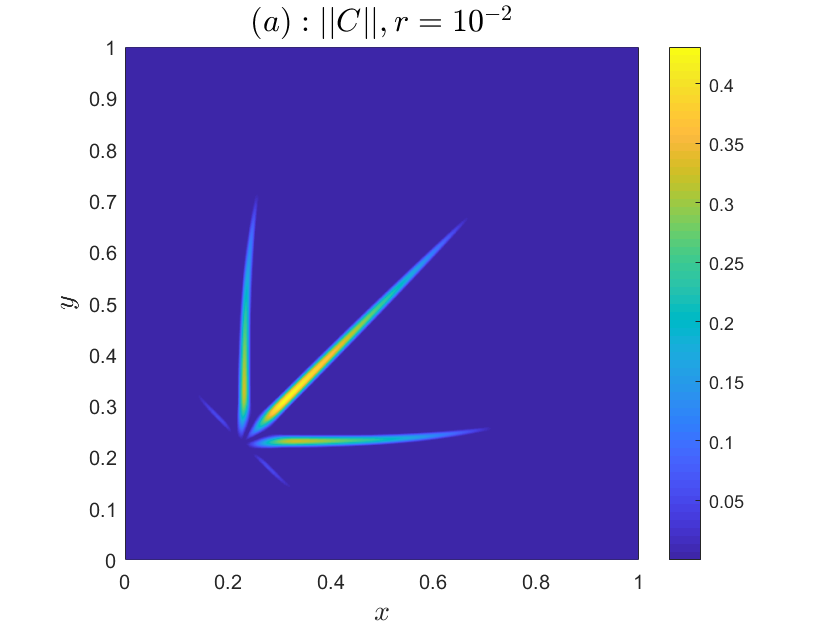}
	\end{minipage}
	\begin{minipage}{.32\textwidth}    \includegraphics[width=1.1\textwidth]{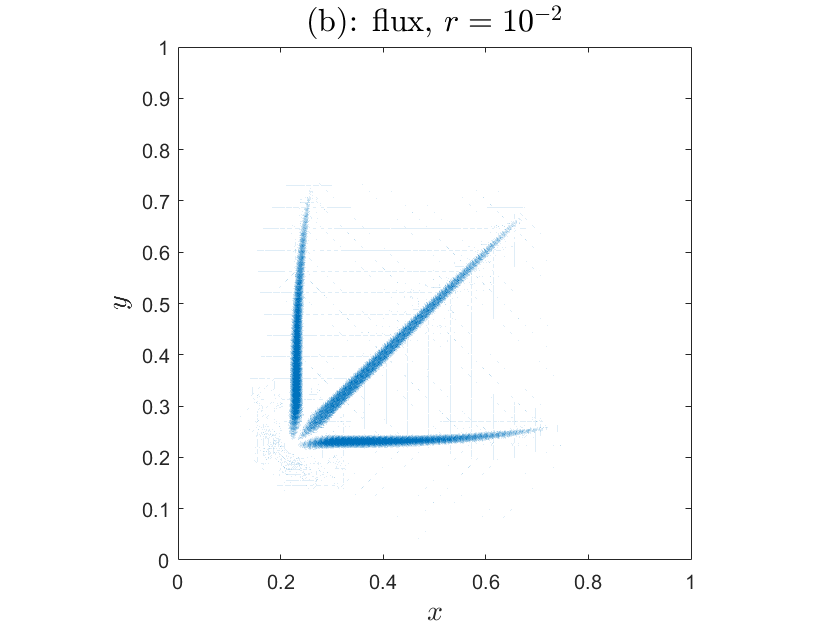}
	\end{minipage}
	\begin{minipage}{.32\textwidth}    \includegraphics[width=1.1\textwidth]{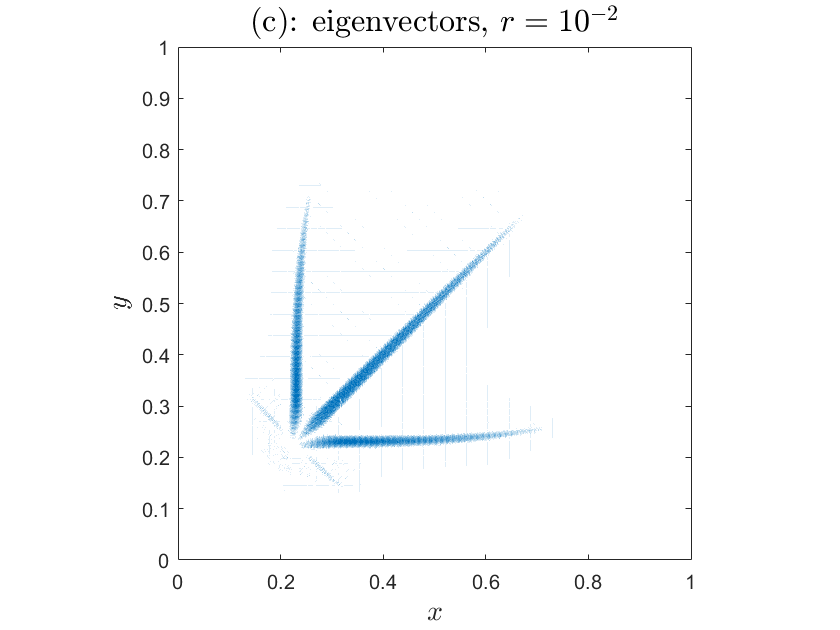}
	\end{minipage}
	\begin{minipage}{.32\textwidth}    \includegraphics[width=1.1\textwidth]{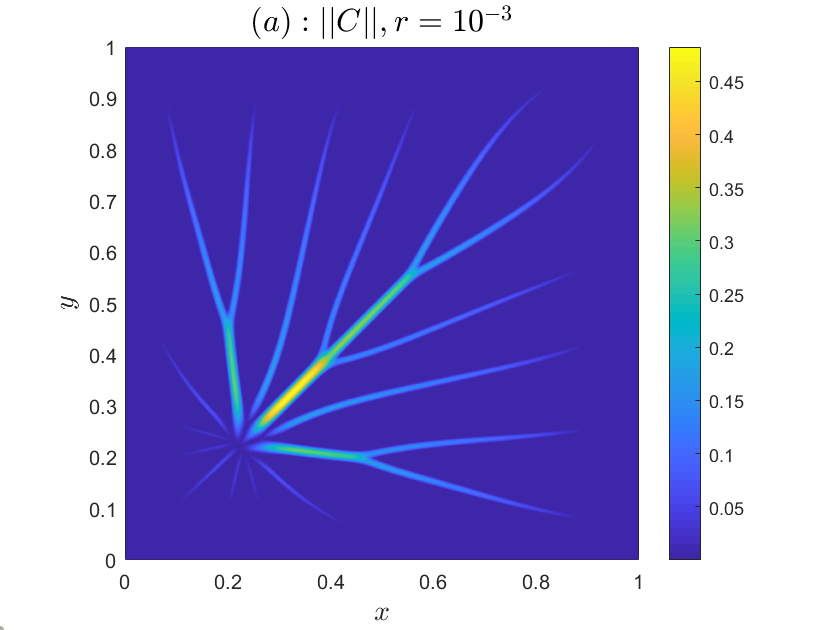}
	\end{minipage}
	\begin{minipage}{.32\textwidth}    \includegraphics[width=1.1\textwidth]{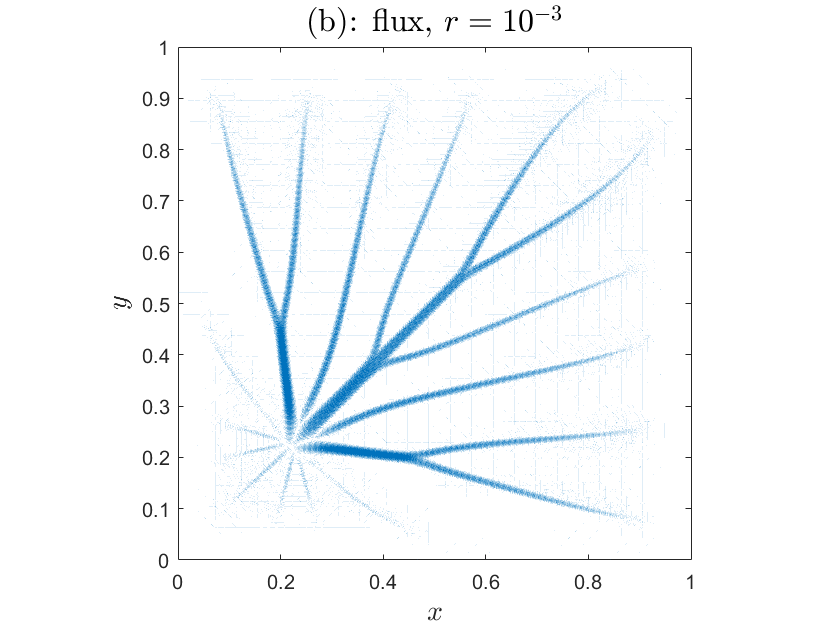}
	\end{minipage}
	\begin{minipage}{.32\textwidth}    \includegraphics[width=1.1\textwidth]{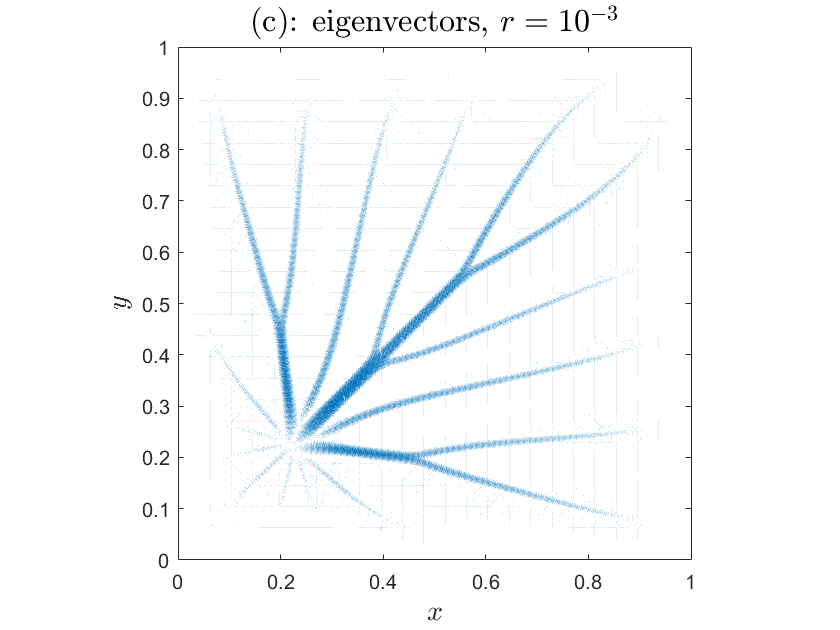}
	\end{minipage}
	\begin{minipage}{.32\textwidth}    \includegraphics[width=1.1\textwidth]{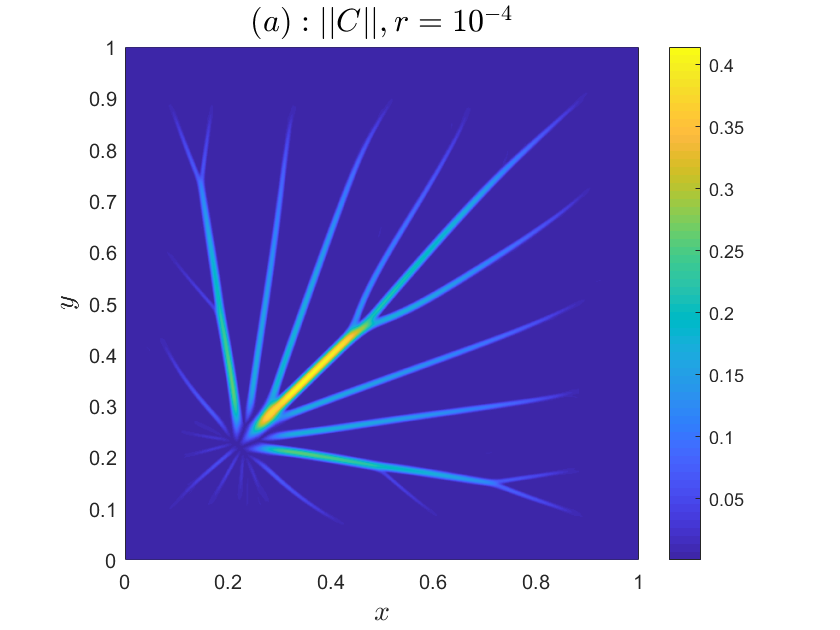}
	\end{minipage}
	\begin{minipage}{.32\textwidth}    \includegraphics[width=1.1\textwidth]{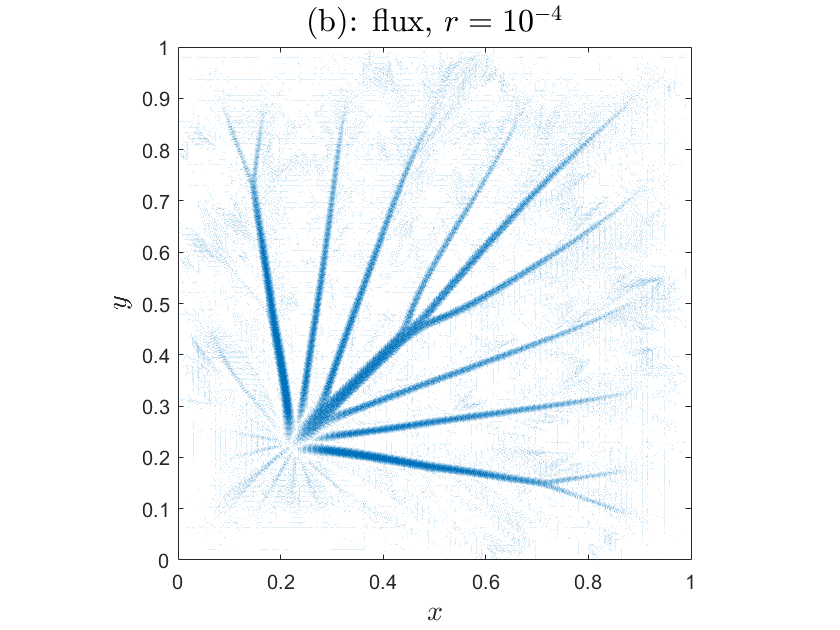}
	\end{minipage}
	\begin{minipage}{.32\textwidth}    \includegraphics[width=1.1\textwidth]{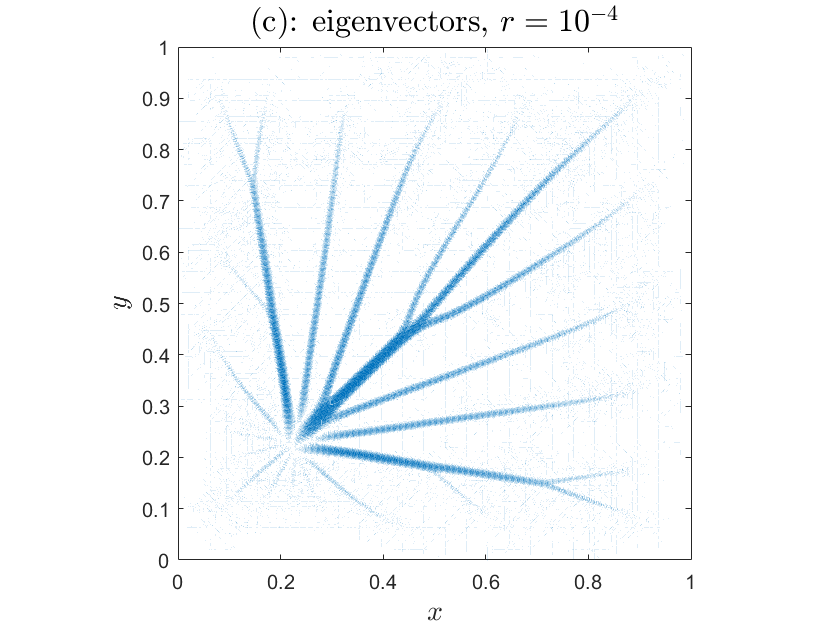}
	\end{minipage}    
	\caption{\textit{In this figure we show three different quantities of the same computations, with the parameters defined in Eq.~\eqref{eq_ic1_def}:  the { Frobenius norm of the conductivity} at final time (a), {the two components of the flux
				{$||\mathbb{C}\nabla p||$}} at final time (b) and the eigenvectors associated to the greatest eigenvalues { of the tensor $\mathbb{C}$ } in absolute value (c). The first row is for $r = 10^{-2}$, the second row for $r = 10^{-3}$ and the last one for $r = 10^{-4}.$ The rest of the data are defined in Eqs.~(\ref{eq_ic1_def}-\ref{eq_choice_param}).}} 
	\label{fig_ic1}
\end{figure}

\begin{figure}[H]
	\centering
	\begin{minipage}{.32\textwidth}    \includegraphics[width=1.1\textwidth]{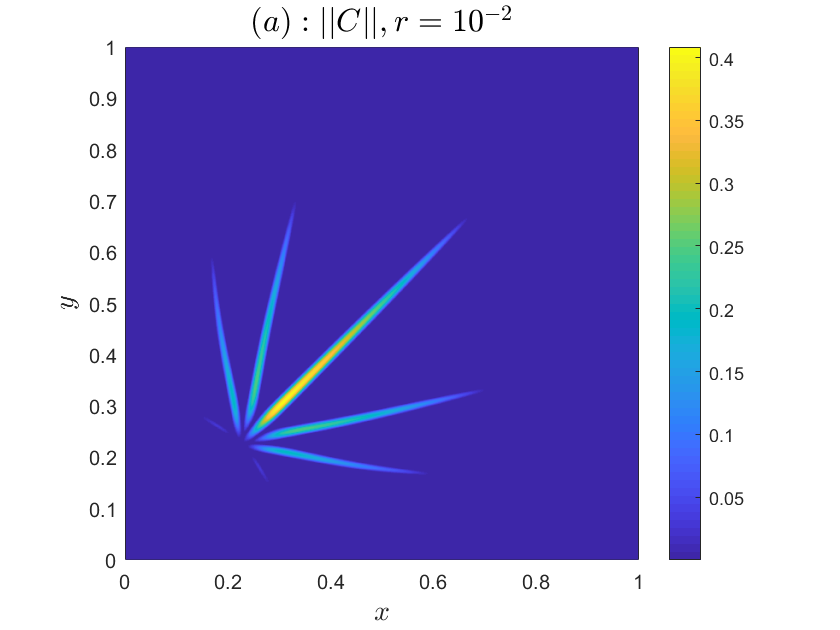}
	\end{minipage}
	\begin{minipage}{.32\textwidth}    \includegraphics[width=1.1\textwidth]{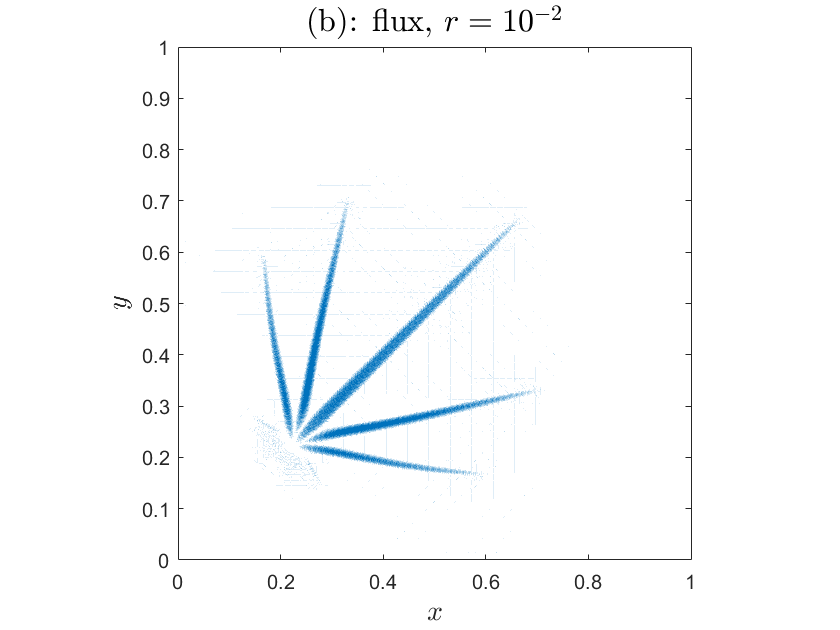}
	\end{minipage}
	\begin{minipage}{.32\textwidth}    \includegraphics[width=1.1\textwidth]{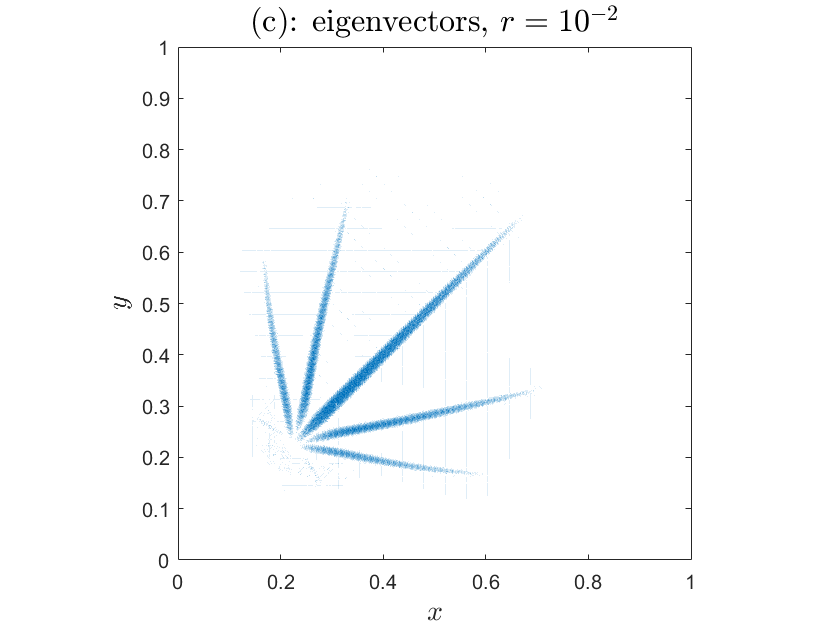}
	\end{minipage}
	\begin{minipage}{.32\textwidth}    \includegraphics[width=1.1\textwidth]{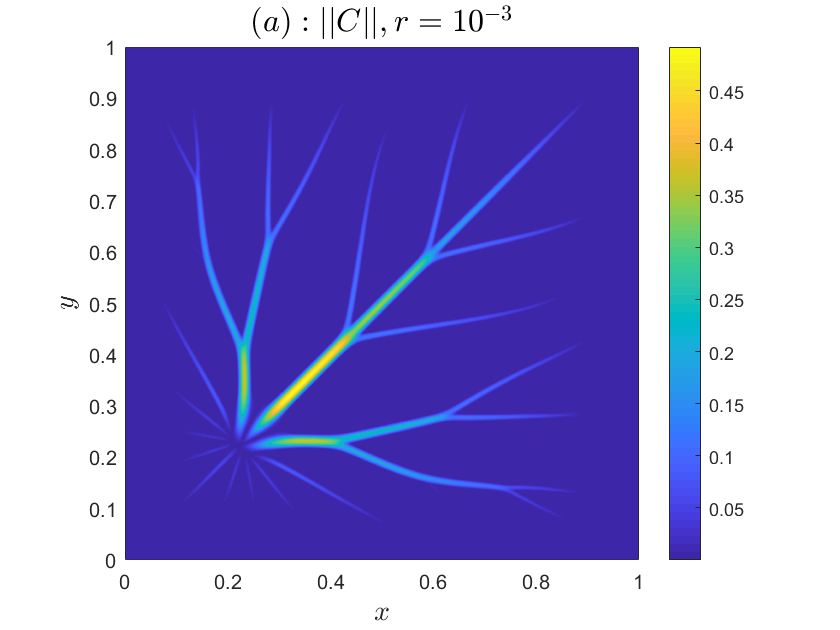}
	\end{minipage}
	\begin{minipage}{.32\textwidth}    \includegraphics[width=1.1\textwidth]{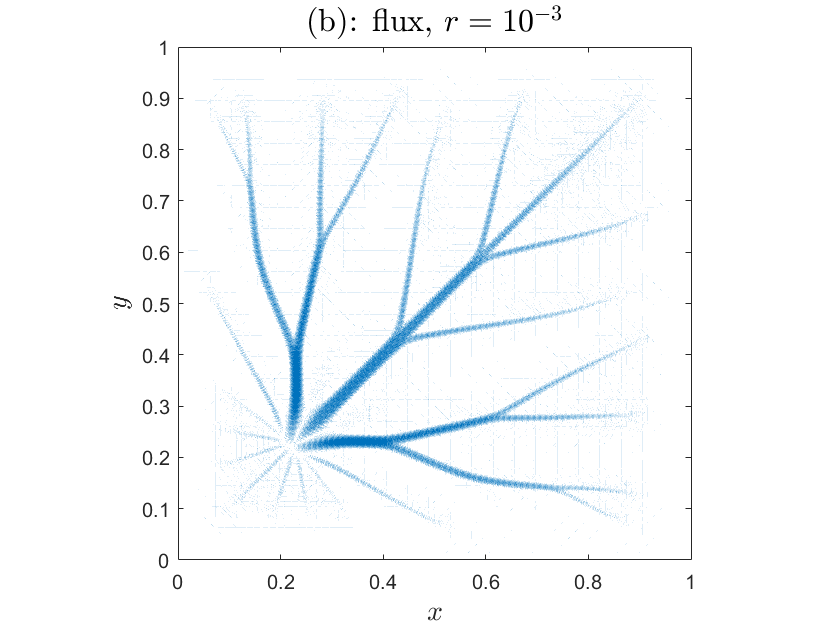}
	\end{minipage}
	\begin{minipage}{.32\textwidth}    \includegraphics[width=1.1\textwidth]{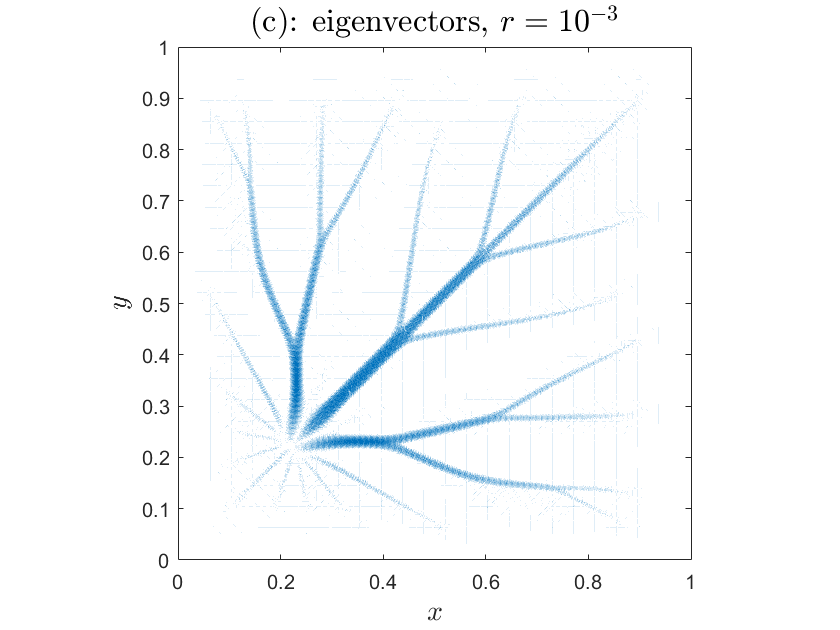}
	\end{minipage}
	\begin{minipage}{.32\textwidth}    \includegraphics[width=1.1\textwidth]{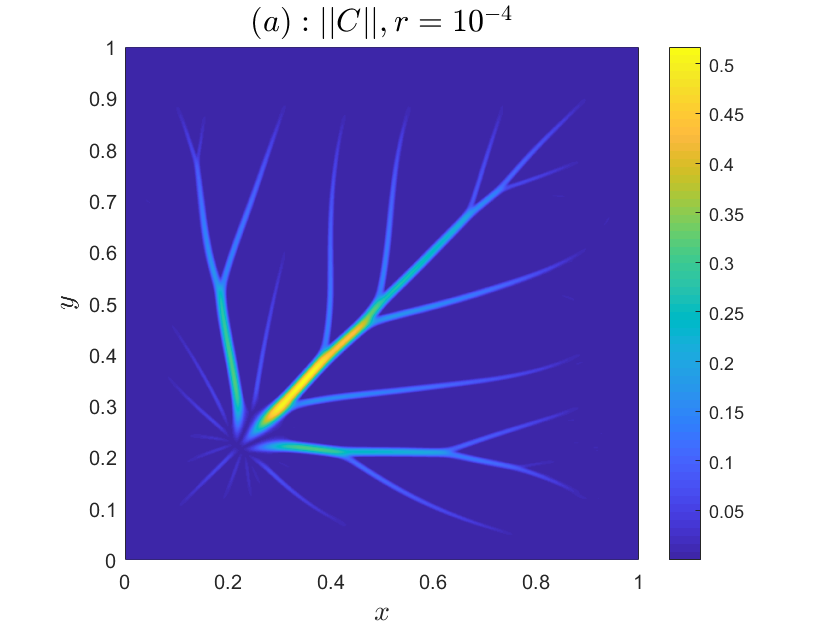}
	\end{minipage}
	\begin{minipage}{.32\textwidth}    \includegraphics[width=1.1\textwidth]{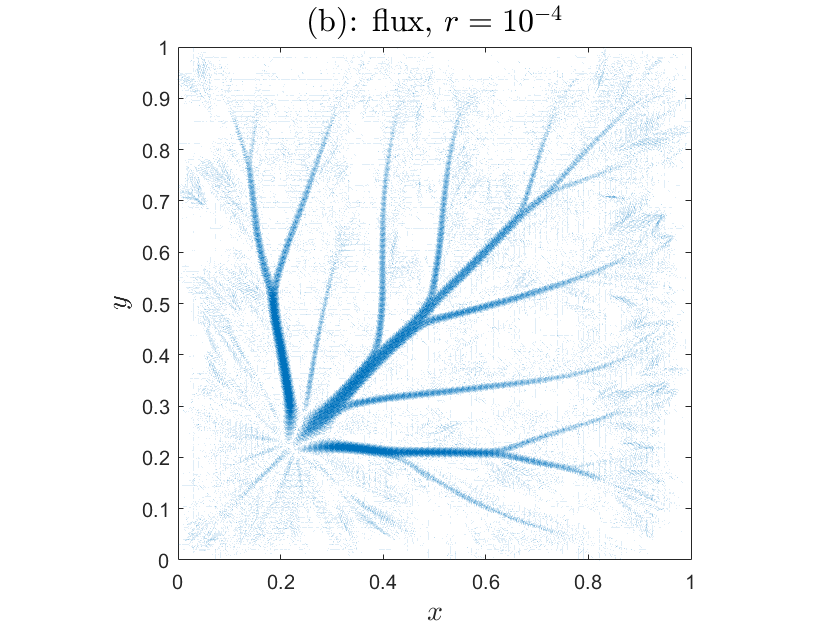}
	\end{minipage}
	\begin{minipage}{.32\textwidth}    \includegraphics[width=1.1\textwidth]{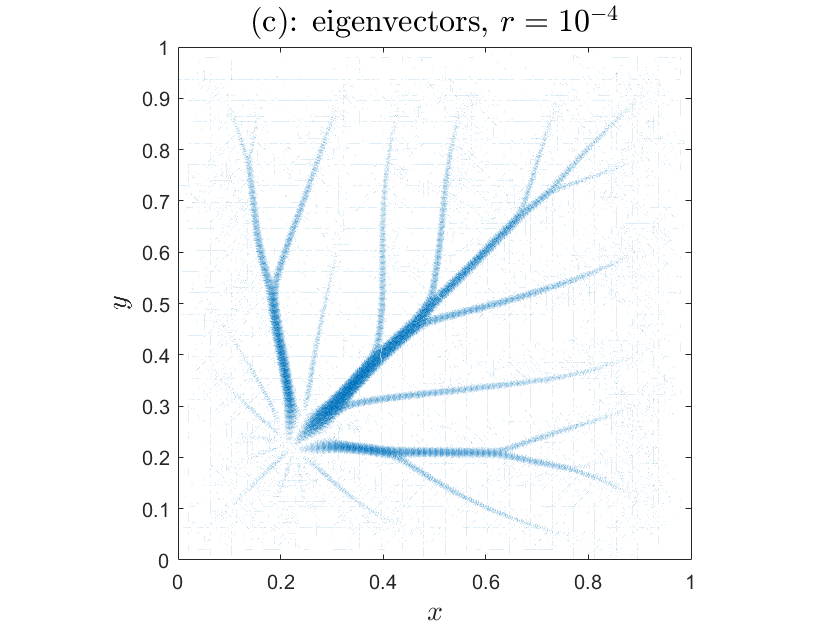}
	\end{minipage}    
	\caption{\textit{In this figure we show three different quantities of the same computations, with the parameters defined in Eq.~\eqref{eq_ic3_def}: the { Frobenius norm of the conductivity} at final time (a), {the two components of the flux
				{$||\mathbb{C}\nabla p||$}} at final time (b) and the eigenvectors associated to the greatest eigenvalue { of the tensor $\mathbb{C}$ } in absolute value (c). The first row is for $r = 10^{-2}$, the second row for $r = 10^{-3}$ and the last one for $r = 10^{-4}.$ The rest of the data are defined in Eqs.~(\ref{eq_ic1_def}-\ref{eq_choice_param}).}} 
	\label{fig_ic3}
\end{figure}

\begin{figure}[h]
	\centering
	\begin{minipage}{.48\textwidth}\includegraphics[width=1\textwidth]{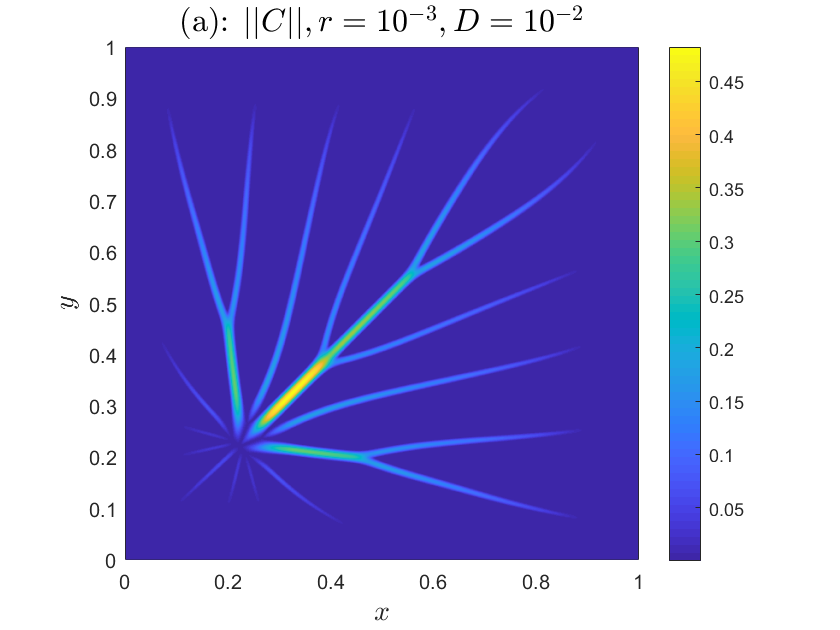}
	\end{minipage}
	\begin{minipage}{.48\textwidth}   \includegraphics[width=1\textwidth]{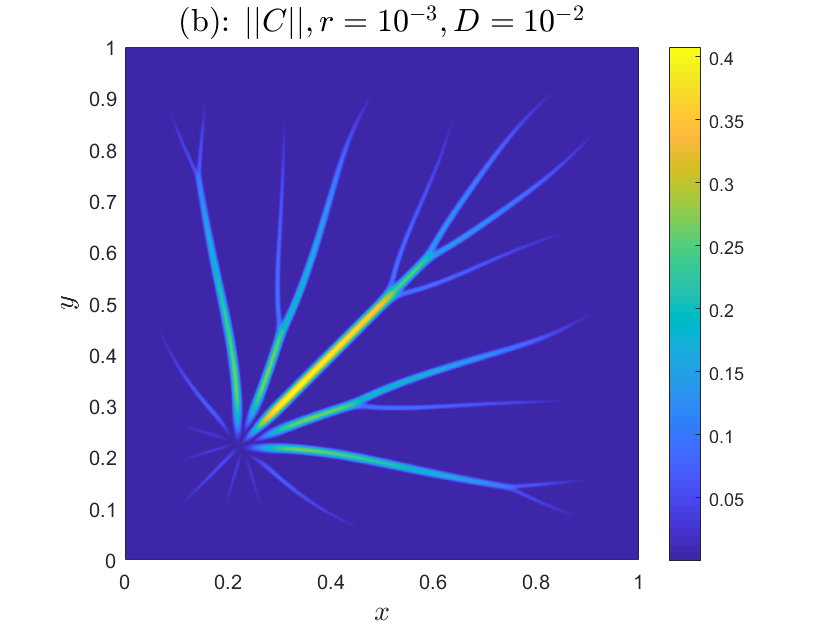}
	\end{minipage}
	\begin{minipage}{.48\textwidth}\includegraphics[width=1\textwidth]{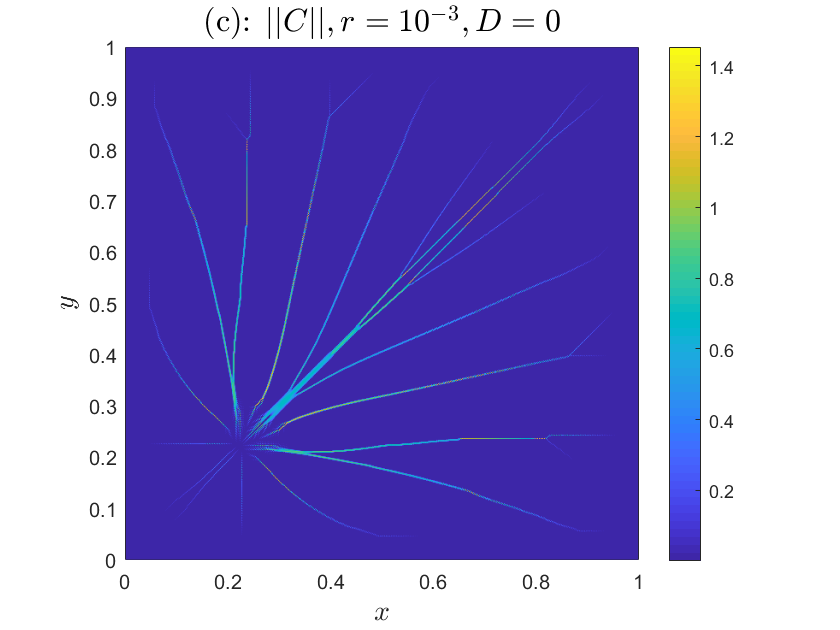}
	\end{minipage}
	\begin{minipage}{.48\textwidth}   \includegraphics[width=1\textwidth]{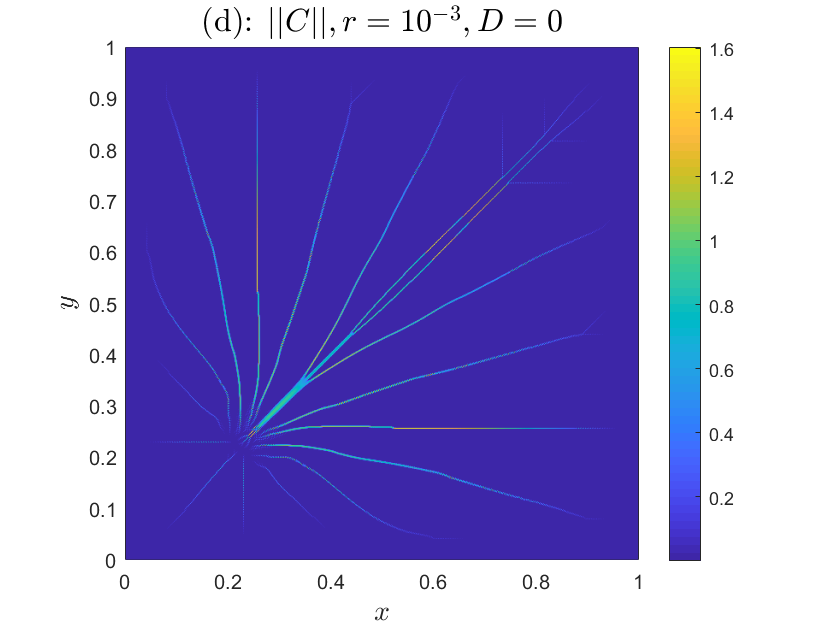}
	\end{minipage}    
	\caption{\textit{Comparison between the module of the solutions, at final time, with different initial conditions and diffusion coefficients $D$. In (a) and (c) the initial condition is defined in Eq.~\eqref{eq_ic1_def}, while in (b) and (d) we have zero initial condition defined in Eq.~\eqref{eq_ic2_def}. In plots (a) and (b) the diffusion coefficient is {$ D = 10^{-2}$}, while in plots (c) and (d) we have zero diffusivity, i.e., $D = 0$.}}   
	\label{fig_comp_ic0}
\end{figure}

To better understand how strong is the dependency, in Figs. \ref{fig_ic1}-\ref{fig_ic3} we show the plots of three different features of the solutions, the {Frobenius norm of the tensor $\mathbb{C}$} at final time (a), {the two components of the flux {$||\mathbb{C}\nabla p||$}} at final time (b) and the eigenvectors associated to the greatest eigenvalues in absolute value {of the tensor $\mathbb{C}$ (i.e., since we have discrete values of the four components of the tensor $C^{(k,l)}_{i,j}, \, k,l = 1,2, \, (i,j) \in \{1,\cdots,N\}^2$, for all $(i,j)$ we consider the 
	$(i-th,j-th)$ term of the tensor $\mathbb{C}_{i,j} = (C^{(1,1)}_{i,j},C^{(1,2)}_{i,j};C^{(2,1)}_{i,j},C^{(2,2)}_{i,j})$, and the corresponding couple of eigenvalues $(\textsc{e}_1,\textsc{e}_2)$ and eigenvectors $(\textsc{v}_1,\textsc{v}_2)$. In the plot we show the two components of the normalized eigenvectors associated to the greatest eigenvalues in absolute value)} (c), for $r = 10^{-2}, 10^{-3}, 10^{-4}$. We observe that the symmetry of the solution breaks for $r = 10^{-4}$, as we expected from Fig.~\ref{fig_condition}, in which ${\rm asymm}(\mathbb{C}) \approx 10^{-1}$.
In Fig.~\ref{fig_comp_ic0} we show the comparison between the solutions obtained with two different initial conditions, and choosing $r = 10^{-3}$: Eq.~\eqref{eq_ic1_def} for plots (a) and (c), and Eq.~\eqref{eq_ic2_def} for (b) and (d). Another difference among the plots is the diffusivity. We choose $D = 10^{-2}$ for plots (a) and (b), and $D = 0$ for plots (c) and (d). The main difference between (a) and (b) panels is the number of branches. Choosing zero initial condition, {the solution is more detailed because it has more ramifications} respect to the one with initial condition constant, equal to one.  This is an effect of the permeability tensor field $\mathbb{P}[\mathbb{C}] = \mathbb{C} + r\mathbb{I}$ in the elliptic operator. In panel (a), it is true that $ \mathbb{C} > r = 10^{-3}$ at initial times, thus, even if we decrease $r$, the $\mathbb{C}$ component is the one prevailing. On the contrary, in panel (b), choosing the initial condition equal to 0, {the parameter $r$ is the one that is} leading the elliptic operator at initial times, and being smaller than the module of the other initial condition, we are able to see more features. Regarding the case of zero diffusivity, we calculated the condition number of the elliptic operators, for the (c) panel ${\rm cond}(L) = 7.030\times 10^{7}$ and for the (d) panel ${\rm cond}(L) = 9.845\times 10^{7}$. This means that, with zero diffusivity, the stiffness of the system increases, and the solutions are less symmetric and accurate.    

\subsection{Accuracy tests: Wasserstein distance}
\label{section_wasserstein}
{
	As can be observed in Fig. \ref{fig_ic1},
	the numerical solution $\mathbb{C}$ is close to zero in a significant portion of the computational domain, away from the branches of the emergent network. Since the shape, position and number of the branches change as the grid is being refined, using an $L^p$ norm is not advantageous in studying the convergence analysis.
	Instead, we propose to use the Wasserstein distance \cite{villani2003topics, villani2009optimal, villani2009wasserstein}, which measures the displacement of the solutions obtained with different numerical grids. For the convenience of the reader, we shall shortly elaborate on how the Wasserstein distance is defined and illustrate its concept with a simple example.}

\begin{figure}[H]
	\centering
	\begin{minipage}{.49\textwidth}  \begin{overpic}[abs,width=\textwidth,unit=1mm,scale=.25]{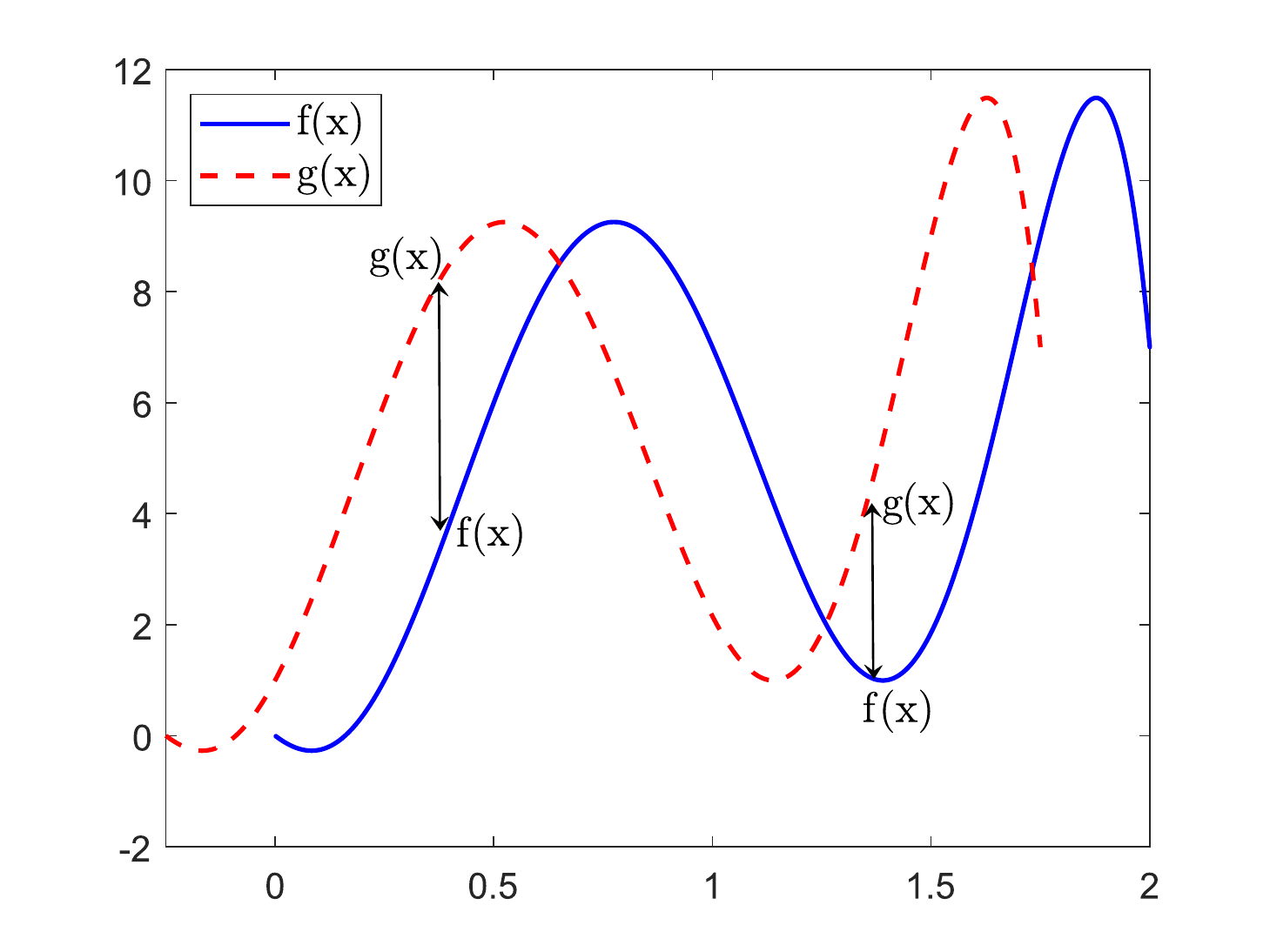}
			\put(17.,47){'vertical distance'}\end{overpic}	
	\end{minipage}
	\begin{minipage}{.49\textwidth} \begin{overpic}[abs,width=\textwidth,unit=1mm,scale=.25]{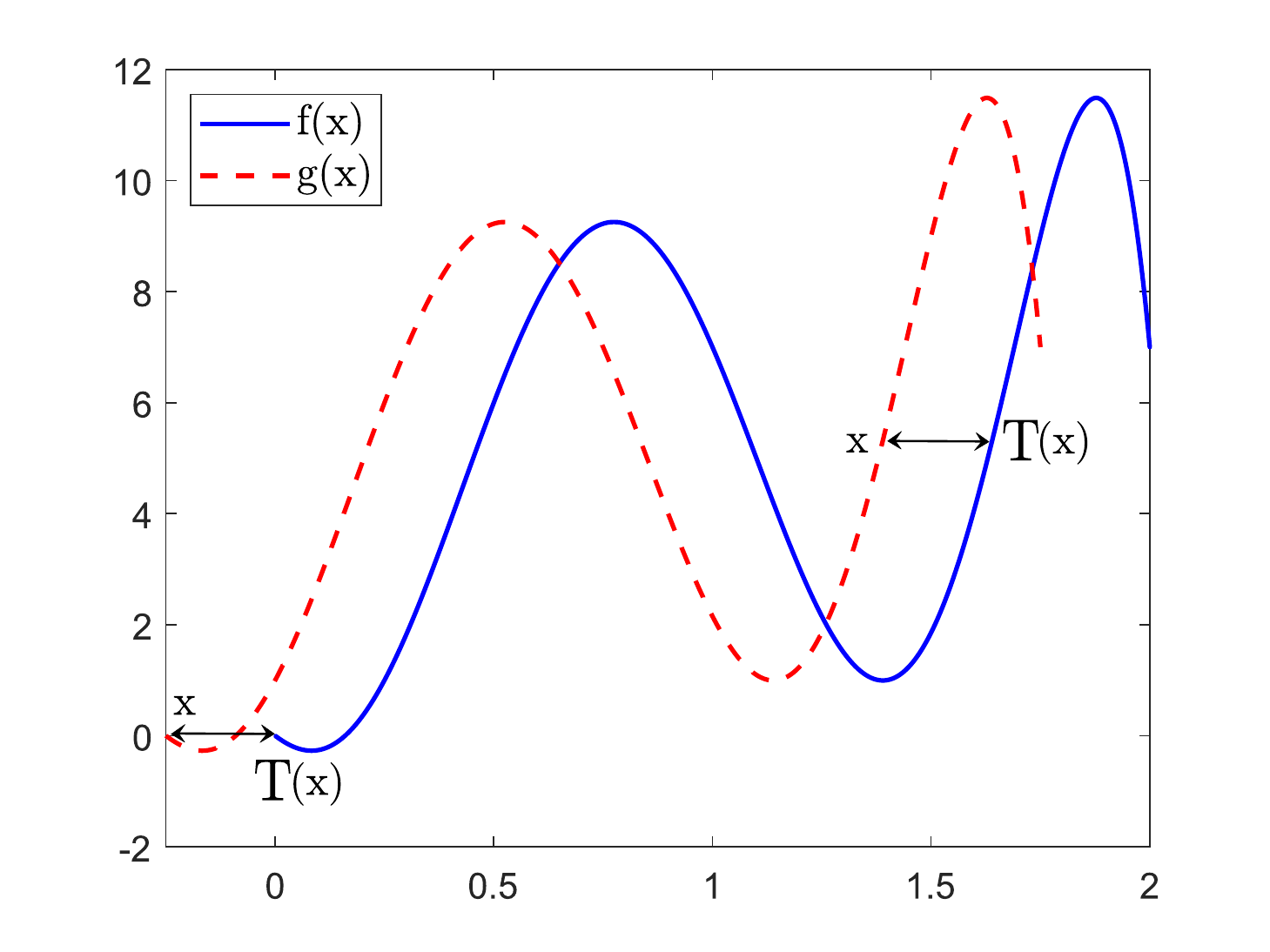}
			\put(17,47){'horizontal distance'}
		\end{overpic}	 
	\end{minipage}
	\caption{{\textit{'Vertical' vs 'horizontal' distances between two generic functions: $f(x), x\in \mathbb{R}$ and its translation $g(x)$, s.t. $g(x) = f(T(x))$ anf and $T(x) = x + {\rm k}, \,{\rm k}\in \mathbb{R}^+$. If we associate the Wasserstein distance to the 'horizontal' distance $'T(x) - x'$, the figure demonstrates that, loosely speaking, it depends more on the displacement of the function than its shape (or 'vertical distance') \cite{santambrogio2017euclidean}.} }}
	\label{fig_WD}
\end{figure}

The Wasserstein metric comes from the idea of moving a distribution of mass, minimizing the average of displacement, see Fig.~\ref{fig_WD} and, e.g., \cite{villani2003topics,villani2009optimal,ambrosio2005gradient}. {In Fig.~\ref{fig_WD} we show a generic function $f(x), x\in \mathbb{R}$ and its translation $g(x)$, such that $g(x) = f(T(x))$ and $T(x) = x + {\rm k}, \,{\rm k}\in \mathbb{R}^+$. With this picture, we address to the 'horizontal' distance (right panel) a 'visual' explanation of the advantage of considering the Wasserstein distance between the two functions, contrary to their 'vertical' one (left panel). }



The Wasserstein distance of order $p\in[1,+\infty)$ between the probability measures $\mu,\nu$ from the metric space $(X,d)$ is defined as
\begin{equation}
	W_p(\mu,\nu) = \Big( \inf_{\pi\in\Pi(\mu,\nu)}  \int_{X\times X} d(x,y)^p\, {\rm d}\pi(x,y)  \Big)^{1/p},
\end{equation}
where $\Pi(\mu,\nu)$ is the set of all transference plans between $\mu$ and $\nu$, i.e., measures on the product space $(X,d)^2$ with marginals $\mu$ and $\nu$, respectively.

\begin{figure}[H]
	\centering
	\begin{minipage}{.49\textwidth}    \includegraphics[width=1.\textwidth]{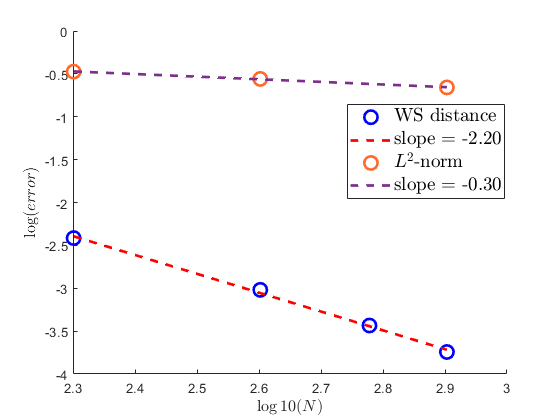}
	\end{minipage}    
	\begin{minipage}{.49\textwidth}
		\includegraphics[width=1\textwidth]{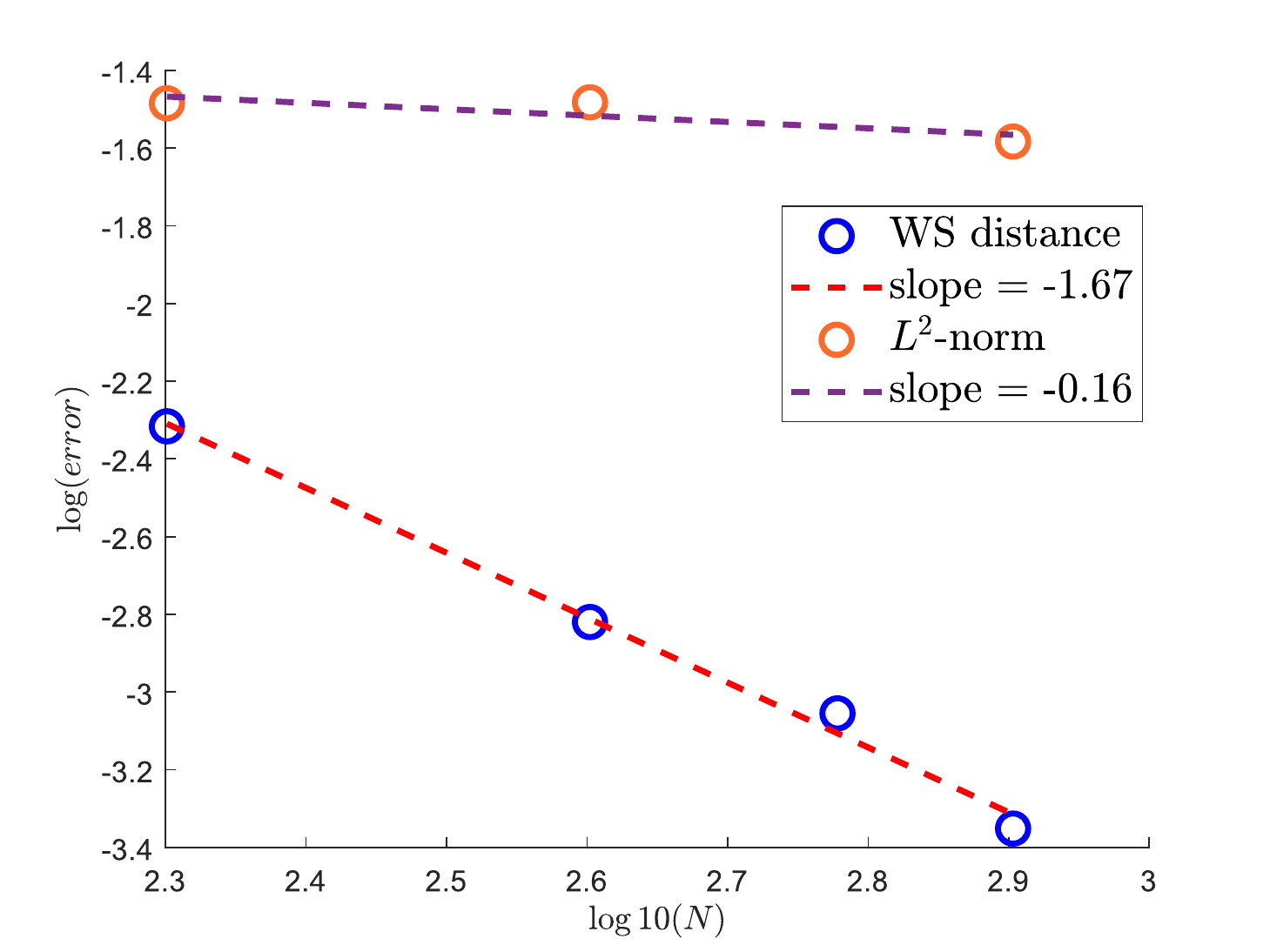}
	\end{minipage}    
	\caption{\textit{Comparison between the accuracy calculated with $L^2$-norm and Wasserstein distance with initial condition defined in Eq.~\eqref{eq_ic1_def}, with $r = 10^{-2}$ (left panel, without extrapolation technique in time) and $r = 10^{-3}$ (right panel, with extrapolation  in time). The rest of the data are defined in Eqs.~\eqref{eq_ic1_def}--\eqref{eq_choice_param}.}}
	\label{fig_acc_comparison}
\end{figure}

In Fig.~\ref{fig_acc_comparison} and Table~\ref{tab_wasserstein} we give a comparison between the accuracy of our numerical simulations calculated with $L^2$-norm and the Wasserstein metric\footnote{{To calculate the Wasserstein distances between two vectors, we make use of the matlab function available at https://github.com/nklb/wasserstein-distance.}}  for $r = 10^{-2}$, on the left panel, and $r = 10^{-3}$, on the right panel.
We used the initial datum \eqref{eq_ic1_def} and the parameter settings specified in \eqref{eq_ic1_def}--\eqref{eq_choice_param}.
We observe an improvement of the order of accuracy using the Wasserstein metric. When considering the $L^2$ norm, the error analysis does not seem to be affected by the number of points used in the space discretization.

In Table~\ref{tab_wasserstein} we present in the first column, the Wasserstein distance, error$_\mathcal{W}$, between two solution with {spatial steps equal to $h$ and $2h$}. In the second column, we calculate the error with the Richardson extrapolation technique, error$_\mathcal{R}$. {It is a technique to estimate the error of a numerical scheme, when the order $p$ is known. 
	If we denote by $u_{\rm exa}$ the exact solution of the system, by $u(h)$ the solution that depends on a discretization parameter $h$, and by $u(0)$ the limit solution obtained as $h\to0$, with the assumption that $u(0)=u_{\rm exa}$, then we have:
	\begin{align*}
		u(h) & =  u(0) + C h^p + o(h^p)\\
		u(qh) & =  u(0) + C (qh)^p + o(h^p).
	\end{align*}
	Subtracting the second relation from the first, one has:
	\[
	u(h)-u(qh) = Ch^p(1-q^p)+o(h^p)
	\]
	from which it follows
	\[
	Ch^p = \frac{u(h)-u(qh)}{1-q^p}+o(h^p)
	= u(h)-u(0) + o(h^p)
	\]
	For $p=2$ and $q=1/2$, one has $u(h)-u(0) = \frac{4}{3}(u(h)-u(h/2)) + o(h^2).$
	At this point, we define by \[{\rm error}_{\mathcal{R}} = \frac{||\mathbb{C}(h) - \mathbb{C}(h/2)||_2}{||\mathbb{C}(h)||_2}\]
	the error calculated with Richardson extrapolation, where $||\cdot||_2$ is the $L^2$-norm.}

We observe {a} better rate of convergence when using the Wasserstein distance,
which is related to the fact that the Wasserstein distance depends to a larger extent on the relative displacement of its arguments, i.e., the topological features of the network, rather than on the local values of the solution at each grid point. 



\begin{table}[]
	\caption{\textit{Accuracy tests with Wasserstein distance, ${\rm error}_\mathcal{W}$, and Richardson extrapolation, ${\rm error}_\mathcal{R}$, for two different values of $r$: on the left $r = 10^{-2}$ and on the right $r = 10^{-3}$, at final time $t = 3$. With the Wasserstein distance, we are able to show the second order of the scheme, obtaining a result that is drastically better than the ones obtained with the usual norms. Here, again, we see that for a smaller values of $r$ (right panel), the accuracy of the method gets worse. The rest of the data are defined in Eqs.~(\ref{eq_ic1_def}-\ref{eq_choice_param}).}} 
		\centering
		\begin{tabular}{ccc||ccc}
			\hline
			$N$  &  error$_\mathcal{W}$  &  error$_\mathcal{R}$ & $N$  &  error$_\mathcal{W}$  &  error$_\mathcal{R}$ \\
			\hline
			100  &  -  & - & 100  &  -  & - \\
			\hline
			200 & 3.846e-03  &  3.358e-01 & 200 & 4.827e-03  &  3.273e-02\\
			\hline
			400 & 9.579e-04 &  2.763e-01 & 400 & 1.514e-03 &  3.292e-02\\
			\hline
			800  &  1.8060e-04 &  2.199e-01 & 800  & 4.455e-04 &  2.609e-02\\
			\hline
		\end{tabular} 	
	\label{tab_wasserstein}
\end{table}

\section{Conclusions}
\label{section_conclusion}
In this paper we explored the effects of the background permeability parameter $r$ in the elliptic-parabolic model \eqref{eq_darcy_p_fin}--\eqref{eq_reaction_diff_fin} which represents a PDE framework describing the formation of biological networks. We showed that $r$ is the parameter that influences the condition number of the elliptic operator in \eqref{eq_darcy_p_fin}. When the condition number becomes {very large}, the {asymmetry of the numerical solution increases}.  
As all the parameters of the system, the initial datum and the source function are symmetric, using a symmetric numerical scheme, we expect the solution to be symmetric at each time step. We make use of finite differences scheme to compute the solution of the system, with central differences for the space discretization and a \textit{symmetric-}ADI method in time. Because of the non linearity in the Poisson equation, we make use of a time extrapolation to improve the accuracy convergence.

Moreover, we demonstrate that using the Wasserstein distance to measure the order of convergence of the numerical scheme provides better results than the $L^2$ norm, in particular when considering small values of the background permeability $r$. In \cite{astuto2022comparison} we were able to show the (expected) second order accuracy only in the case of $r = 0.1$. Here, using the Wasserstein distance to measure the error, we see a significant improvement even for much smaller values of the background permeability, $r = 10^{-2}$ and $r = 10^{-3}$.

In a future work we shall improve the numerical method by implementing a monolithic algorithm to solve the system, {with a linearization for the conductivity variable in the Poisson equation for the pressure}. In that way, the two unknowns, conductivity tensor $\mathbb{C}$ and pressure $p$, are solved implicitly at each time step, with the metabolic term that can be easily linearized, {as we already did in this paper}. An application of IMEX schemes will ensure a higher order of accuracy. The monolithic scheme would also improve the stability in time, allowing us to consider larger time steps.

\section*{Declarations}

\subsection*{Funding and Conflict of interest}
The authors declare no external funding and no conflict of interest.

\subsection*{Authors' contribution}
P.M. and J.H. modelled the systems that govern this work. C.A., D.B. and G.R. set the methodologies and numerical schemes that were used. C.A. validated the results and wrote the original draft. C.A., J.H. and G.R. wrote, reviewed, and edited the manuscript. D.B, acquired funding and provided resources. 

\subsection*{Consent to participate and for publication}
All the authors consent their participation and agree to the published version of the manuscript.


\begin{thebibliography}{99}
	\bibitem[1]{malinowski2013understanding} Malinowski, R.: Understanding of leaf development -- the science of complexity. Plants, 2013, Multidisciplinary Digital Publishing Institute. 
	
	\bibitem[2]{sedmera2011function} Sedmera, D.: Function and form in the developing cardiovascular system. Cardiovascular research, 2011,
	Oxford University Press. 
	
	\bibitem[3]{eichmann2005guidance} Eichmann, A., Le Noble F., Autiero  M., Carmeliet P.: Guidance of vascular and neural network formation.  Current opinion in neurobiology, 2005, Elsevier. 
	
	\bibitem[4]{michel1995morphogenesis} Michel O., Biondi J.: Morphogenesis of neural networks.
	Neural Processing Letters, 1995, Springer. 
	
	\bibitem[5]{pohl1986crucial} Crucial role of endothelium in the vasodilator response to increased flow in vivo,
	Pohl, U.,  Holtz, J., Busse R.,  Bassenge, E.: Hypertension. 1986, American Heart Association. 
	
	\bibitem[6]{hacking1996shear}
	Hacking, W.J., VanBavel, E., Spaan, J.A.: Shear stress is not sufficient to control growth of vascular networks: a model study. American Journal of Physiology-Heart and Circulatory Physiology, 1996, American Physiological Society Bethesda, MD. 
	
	\bibitem[7]{pries1998structural}
	Pries, A.R.,  Secomb T.W., Gaehtgens, P.: Structural adaptation and stability of microvascular networks: theory and simulations.  American Journal of Physiology-Heart and Circulatory Physiology, 1998, American Physiological Society Bethesda, MD. 
	
	\bibitem[8]{chen2012haemodynamics}
	Chen, Q.,  Jiang, L., Li, C., Hu, D., Bu, J., Cai D., Du, J.: Haemodynamics-driven developmental pruning of brain vasculature in zebrafish. 2012, Public Library of Science San Francisco, USA. 
	
	\bibitem[9]{hu2012blood}
	Hu, D., Cai D.,  Rangan, A.V.: Blood vessel adaptation with fluctuations in capillary flow distribution. PLoS ONE, 2012, Public Library of Science San Francisco, USA. 
	
	\bibitem[10]{hu2013adaptation} Hu, D.,  Cai, D.: Adaptation and optimization of biological transport networks. Physical review letters, 2013, APS. 
	
	\bibitem[11]{hu2013optimization}
	Hu D.: Optimization, adaptation, and initialization of biological transport networks.   Notes from lecture, 2013
	
	\bibitem[12]{marko_perthame} {Fang, D.,  Jin, S., Markowich, P., Perthame, B.:} Implicit and {Semi-implicit} {Numerical} {Schemes} for the {Gradient} {Flow} of the {Formation} of {Biological} {Transport} {Networks}. {The SMAI journal of computational mathematics}, {2019}. 
	
	\bibitem[13]{marko_perthame_2}
	Haskovec, J., Markowich P., Perthame, B.: Mathematical Analysis of a PDE System for Biological Network Formation. Communications in Partial Differential Equations, 2015, Taylor \& Francis. 
	
	\bibitem[14]{marko_perthame_schlo}
	Haskovec, J., Markowich, P.,  Perthame, B., Schlottbom, M.: Notes on a PDE system for biological network formation. Nonlinear Analysis, 2016. 
	
	\bibitem[15]{marko_albi}
	Albi, G., Artina, M.,  Foransier, M., Markowich, P.: Biological transportation networks: Modeling and simulation. Analysis and Applications, 2016. 
	
	\bibitem[16]{albi_burger}
	Albi, G.,  Burger, M., Haskovec, J.,  Markowich,  P.,  Schlottbom, M.: Continuum Modelling of Biological Network Formation. Birkh\"auser-Springer (Boston), 2017, Active Particles Vol.I - Theory, Models, Applications
	
	\bibitem[17]{portaro}
	Haskovec, J., Markowich, P., Portaro, S.: Emergence of biological transportation networks as a self-regulated process. Discrete and Continuous dynamical system, 2022, American Institute of Mathematical Sciences (AIMS). 
	
	\bibitem[18]{hu2019optimization}
	Hu D., Cai, D.: An optimization principle for initiation and adaptation of biological transport networks. Communications in Mathematical Sciences, 2019, International Press of Boston
	
	
	\bibitem[19]{marko_pilli}
	Haskovec, J., Markowich P., Pilli, G.: Tensor PDE model of biological network formation.  Communications in Mathematical Sciences, 2022, International Press of Boston
	
	
	\bibitem[20]{astuto2022comparison}
	Astuto, C., Boffi, D., Haskovec, J., Markowich, P., Russo, G.: Comparison of Two Aspects of a PDE Model for Biological Network Formation. Mathematical and Computational Applications, 2022. 
	
	\bibitem[21]{carrillo}
	Carrillo, J.A., Toscani, G.: Wasserstein metric and large-time asymptotics of nonlinear diffusion equations. New Trends in Mathematical Physics, 2004. 
	
	\bibitem[22]{otto1996double}
	Otto, F.: Double degenerate diffusion equations as steepest descent. 1996, https://books.google.com.sa/books?id=oxLdGwAACAAJ, Bonn University
	
	\bibitem[23]{ottoF}
	Otto, F.: The geometry of dissipative evolution equations: the porous medium equation. Communications in Partial Differential Equations, 2001, Taylor \& Francis. 
	
	
	\bibitem[24]{astuto2023finite}
	Astuto, C., Boffi, D., Credali, F.: Finite element discretization of a biological network formation system: a preliminary study. arXiv preprint, 2023. 
	
	\bibitem[25]{alcubierre1992time}
	Alcubierre, M., Schutz, B.F.: Time Symmetric Adi and Casual Reconnection. International Workshop on Numerical Relativity, 1992, Cambridge University Press
	
	
	\bibitem[26]{peaceman_rachford}
	Peaceman, D.W.,  Rachford, H.H.Jr.: The numerical solution of parabolic and elliptic equations. Journal of the Society for industrial and Applied Mathematics, 1955
	
	
	\bibitem[27]{CiCP-31-707}
	Raudino, A., Grassi, A., Lombardo, G., Russo, G., Astuto, C., Corti, M.: Anomalous Sorption Kinetics of Self-Interacting Particles by a Spherical Trap.  Communications in Computational Physics, 2022
	
	
	\bibitem[28]{villani2003topics}
	Villani, C.: Topics in optimal transportation.(books). OR/MS Today, 2003, Institute for Operations Research and the Management Sciences
	
	\bibitem[29]{villani2009optimal}
	Villani, C.: Optimal transport: old and new, 2009,
	Springer
	
	\bibitem[30]{villani2009wasserstein}
	Villani, C.: The wasserstein distances. Optimal transport, 2009, Springer. 
	
	\bibitem[31]{ambrosio2005gradient}
	Ambrosio, L., Gigli N., Savar{\'e}, G.: Gradient flows: in metric spaces and in the space of probability measures. 2005, Springer Science \& Business Media. 
	
	
	\bibitem[32]{santambrogio2017euclidean}
	Santambrogio, F.: Euclidean, metric, and Wasserstein gradient flows: an overview. Bulletin of Mathematical Sciences, 2017, Springer. 
	
\end{thebibliography}

\end{document}